\documentclass[11pt,reqno,oneside]{amsart}
\usepackage{verbatim,amsmath,amssymb,cite,epsfig,xspace,graphicx,multicol,miller}
\usepackage{srcltx,subfigure}

\numberwithin{equation}{section}

\theoremstyle{plain}

\theoremstyle{definition}

\theoremstyle{remark}
\newtheorem{remark}{Remark}[section]

\DeclareFontFamily{OT1}{pzc}{}
\DeclareFontShape{OT1}{pzc}{m}{it}{<-> s * [1.10] pzcmi7t}{}
\DeclareMathAlphabet{\mathpzc}{OT1}{pzc}{m}{it}
\DeclareMathOperator*{\newargmin}{arg\ local\ min}

\newcommand{\Real}{\mathbb R}
\def\R{\mathbb{R}}
\def\N{\mathbb{N}}

\newcommand{\Int}{{\mathbb Z}}

\def\<{\langle}
\def\>{\rangle}

\def\E{\mathcal{E}}
\newcommand{\F}{{\mathcal F}}

\def\A{\mathcal{A}}
\def\C{\mathcal{C}}

\def\Us{\mathcal{U}}

\def\del{\delta}

\renewcommand{\cases}[1]{\left\{ \begin{array}{ll} #1 \end{array} \right.}

\newcommand{\smfrac}[2]{{\textstyle \frac{#1}{#2}}}

\def\half{\smfrac{1}{2}}

\def\fcc{{\mathcal{L}}}
\def\dirv{{\hkl[1 -1 0]}}
\def\planev{{\hkl(1 -1 0)}}
\def\proj{{P}}
\def\qproj{{Q}}
\def\tri{{\mathcal{T}}}
\def\domain{{\Omega}}
\def\pdomain{{\omega}}
\def\periodic{{periodic}}
\def\energya{\mathcal{E}^{a}}
\def\benergya{\E^{a}}
\def\Ub{{\mathcal{U}}}
\def\Ut{\mathcal{V}}
\def\p{{\phi}}
\def\d{D}
\def\dextra{\bar{D}}
\def\N{{\mathcal{M}}}
\def\M{{\mathcal{R}}}
\newcommand{\pder}[2]{\frac{\partial #1}{\partial #2}}

\newcommand{\admis}{{\rm Adm}(\dsp_0)}
\newcommand{\ignore}[1]{}
\def\atomistic{{\A}}
\def\continuum{{\C}}
\def\bondind{\chi_b}
\def\forceqcf{\F^{qcf}}

\def\dsp{{u}}

\begin{document}

\title
[A Computational and Theoretical Investigation of Quasicontinuum Methods]
{A Computational and Theoretical Investigation of
the Accuracy of Quasicontinuum Methods}
\author{Brian Van Koten}
\author{Xingjie Helen Li}
\author{Mitchell Luskin}
\author{Christoph Ortner}

\address{Brian Van Koten\\
School of Mathematics \\
University of Minnesota \\
206 Church Street SE \\
Minneapolis, MN 55455 \\
U.S.A.} \email{vank0068@umn.edu}
\address{Xingjie Helen Li\\
School of Mathematics \\
University of Minnesota \\
206 Church Street SE \\
Minneapolis, MN 55455 \\
U.S.A.} \email{lixxx835@umn.edu}
\address{Mitchell Luskin \\
School of Mathematics \\
University of Minnesota \\
206 Church Street SE \\
Minneapolis, MN 55455 \\
U.S.A.} \email{luskin@umn.edu}

\address{Christoph Ortner\\
Mathematical Institute \\
University of Oxford \\
24-29 St Giles'\\
Oxford OX1 3LB \\
UK} \email {ortner@maths.ox.ac.uk}

\thanks{
This work was supported in part by DMS-0757355,
 DMS-0811039, the PIRE Grant OISE-0967140, the Institute for Mathematics and
Its Applications, and
 the University of Minnesota Supercomputing Institute.
This work was also supported by the Department of Energy under
Award Number {DE-SC}0002085. CO was supported by the EPSRC grant
EP/H003096/1 ``Analysis of Atomistic-to-Continuum Coupling
Methods.'' }



\date{\today}

\begin{abstract}{
We give computational results to study the accuracy
 of several quasicontinuum methods for two benchmark problems
 --- the stability of a Lomer dislocation pair under shear and
 the stability of a lattice to plastic slip under tensile
 loading. We find that our theoretical analysis of the accuracy
 near instabilities for one-dimensional model problems can
 successfully explain most of the computational results for
 these multi-dimensional benchmark problems. However, we also
 observe some clear discrepancies, which suggest the need for
 additional theoretical analysis and benchmark problems to more
 thoroughly understand the accuracy of quasicontinuum methods.
}
\end{abstract}

\maketitle

\section{Introduction}
Multiphysics model coupling has captured the excitement of the engineering
research community for its potential to make possible the numerical simulation of heretofore
computationally inaccessible multiscale problems.  The capability to assess the accuracy
of these multiphysics methods is crucial to both the verification of existing methods
and the development of improved methods.

During the past several years, a theoretical basis has been
developed for estimating the error of atomistic-to-continuum
coupling methods
~\cite{cluster09,gavini10,vankoten.blended,eam.qc,badia:onAtCcouplingbyblending,ortner:qnl1d,prudhomme:modelingerrorArlequin,Dobson:2008c,Dobson:2008b,mingyang,LinP:2003a,LinP:2006a,Ortner:2008a,E:2006,Gunzburger:2008a,LuskinXingjie,Legoll:2005,makridakis10,dobs-qcf2,doblusort:qcf.stab,dobsonluskin08,dobsonluskin07}.
However, the accurate computation of lattice instabilities such
as dislocation formation and movement, crack propagation, and
plastic slip are primary goals of atomistic-to-continuum
coupling methods; and a theoretical analysis of the accuracy of
atomistic-to-coupling methods {\em up to the onset of lattice
instability} of the atomistic energy has thus far only been
achieved for one-dimensional model
problems~\cite{doblusort:qce.stab,doblusort:qcf.stab,ortner:qnl1d,Dobson:2008b,belik10,Ortner:2008a}.
This theoretical analysis and corresponding numerical
experiments have given a precise understanding of the varying
accuracy of these methods for simple model problems,  but it is
not known to what degree these errors are significant for the
multi-dimensional problems of scientific and technological
interest.

We have developed a benchmark test code to study the
atomistic-to-continuum coupling error for a face-centered cubic
(FCC) crystal. Following the benchmark investigation of Miller
and Tadmor~\cite{Miller:2008}, the displacement $U\left(x_1,\,x_2,\,x_3\right)$
of the atoms
from their reference lattice positions is
constrained by the atomistic analogue of continuum ``plane
strain'' symmetry:
\begin{equation*}
U=\left(U_1,\,U_2,\,0\right)\quad\text{and}\quad
U\left(x_1,\,x_2,\,x_3\right)=u\left(x_1,\,x_2\right)
\end{equation*}
where the crystal coordinates are given in terms
of the basis vectors defined using
Miller indices\cite{ashcroft} by $e_1=$\hkl[1 1 0],
$e_2=$\hkl[0 0 1], and $e_3=$\hkl[1 -1 0].  We note that
the computation of the atomistic energy and forces for an atomistic displacement
with plane strain symmetry requires the summation of the interaction
of each atom with neighboring atoms in three-dimensional space.

We have recomputed the benchmark study of Miller and Tadmor
\cite{Miller:2008} for the Lomer dislocation dipole to more
clearly separate the errors due to continuum modeling error,
atomistic-to-continuum coupling error, and solution error. We
allow a more general atomistic domain that is fully surrounded
by a continuum region.  We have also investigated plastic slip
for tensile loading. Further investigation is underway to study
the coarsening error in the continuum region and the far-field
boundary condition error.

We have found that the {\em patch-test consistent}
quasi-nonlocal (QNL) method ~\cite{Shimokawa:2004} gave a more
accurate critical shear strain than the popular ``ghost force
correction'' quasicontinuum coupling method
~\cite{Miller:2003a,Shenoy:1999a} tested in~\cite{Miller:2008}
(see Table~\ref{tab: critical strain for dipole} and
Figure~\ref{fig: error for dipole}). This result confirms our
earlier theoretical analysis of the dependence of the critical
shear strain error on the accuracy of the
atomistic-to-continuum coupling~\cite{doblusort:qce.stab}.

The motivation for the ghost force correction method (equation
\eqref{eqn: gfc iteration} below) was to correct the large coupling
errors of the original quasicontinuum method (QCE) ~\cite{Ortiz:1996a}
by applying a ``dead load'' correction
~\cite{Miller:2003a,Shenoy:1999a,Miller:2008}.  We have recast the
ghost force correction method in a numerical analysis setting as an
iterative method to solve the equilibrium equations for the
force-based quasicontinuum (QCF) method
~\cite{dobsonluskin07,dobsonluskin08}, and we have proven for a
one-dimensional model problem that the ghost force correction
(QCE-QCF) method is inaccurate near lattice instabilities because it
uses the {\em patch-test inconsistent} QCE method as a preconditioner
for the QCF equilibrium equations~\cite{qcf.iterative}. {
The benchmark tests we present here confirm that the QCE-QCF method is
not as accurate as the QNL method near lattice instabilities.

By contrast, we found that the $w^{1,\infty}$ error was smaller for
the QNL method than for the QCE-QCF method far from lattice instabilities
(see Figures~\ref{fig: error for dipole}-\ref{fig: error for
  tension}), which is contrary to expectations based} { on our
analyses of a 1D model} { in
\cite{dobs-qcf2,qcf.iterative,ortner:qnl1d,Dobson:2008b}, where we
have shown that the QCF method (the limit of the QCE-QCF iteration) is
a more accurate approximation than the QNL method.

Another result of our numerical experiments that contradicts our 1D
analysis \cite{qcf.iterative} is that near a slip instability the
QCE-QCF method has comparable accuracy to the QNL-QCF method, which
uses the QNL energy as a preconditioner. To explain this, we note that
our 1D analysis in \cite{qcf.iterative} can be considered a good model
for cleavage fracture, but not for the slip instabilities studied in
the present paper.  We are currently attempting to develop a 2D
benchmark test for cleavage fracture to demonstrate that the QCE-QCF
method may be less accurate than the QNL-QCF method near some types of
lattice instabilities.}


\section{The atomistic and quasicontinuum models}
\label{sec:model}

\subsection{A model for plane strain in the face centered cubic lattice}
\label{sec:planestrain_fcc}

\begin{figure}
\centering
\includegraphics[width = 10cm]{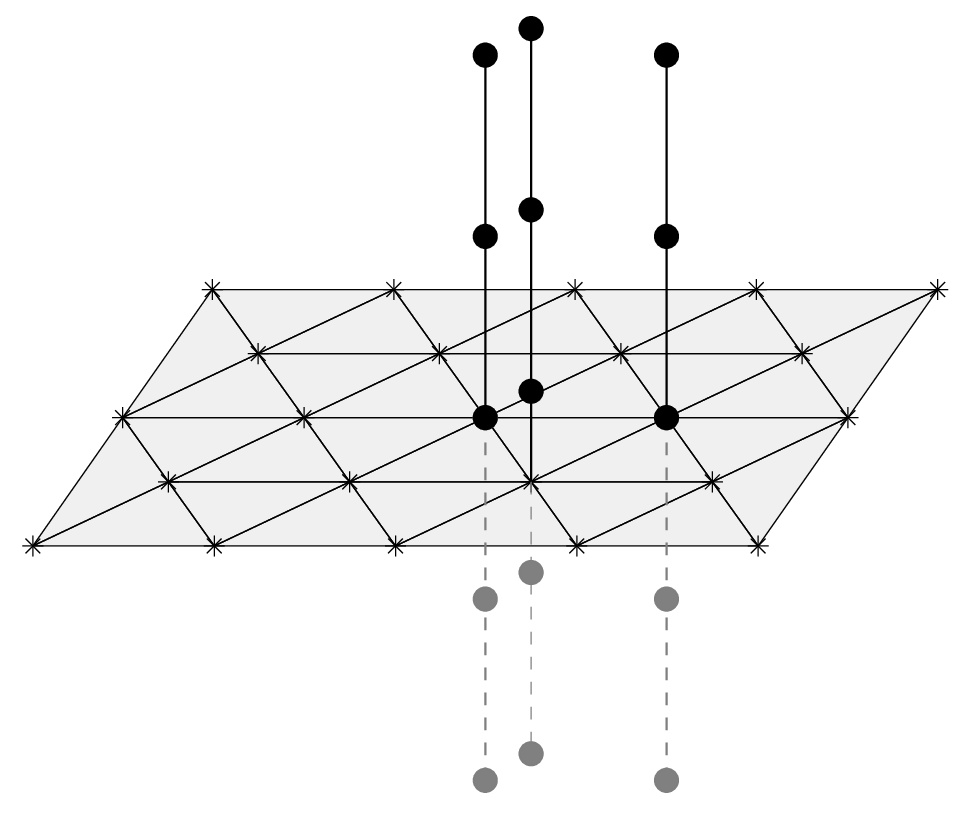}
\caption{Illustration of the FCC lattice and the projection onto the
  $\planev$ plane.  The grey shaded region represents the $\planev$
  plane, the asterisks are the points in the triangular lattice
  $\tri$, and the disks are the atoms in the FCC lattice
  $\mathcal{L}$.  The vertical lines emanating from the points in
  $\tri$ represent the columns of atoms in the FCC lattice which
  project onto points of $\tri$. The reader should imagine that there
  is a vertical column of atoms emanating from each point of $\tri$.}
\label{fig: fcc projection figure}
\end{figure}

Let $\fcc$ be a face centered cubic (fcc) lattice \cite{ashcroft} with cube side length
$a$ and nearest-neighbor distance $a/\sqrt 2$ . Our plane strain model
is most easily derived by viewing the fcc lattice $\fcc$ as being generated by
the primitive lattice vectors
\begin{equation} \label{prim}
a_1=\frac{a}{2} (0,\,1,\, 1),\quad
a_2=\frac{a}{2} (1,\,0,\, 1),\quad\text{and}\quad a_3=\frac{a}{2} (1,\,1,\, 0).
\end{equation}
The cubic supercells are then generated by the cubic lattice vectors
\begin{equation}
\label{cubic}
\begin{gathered}
A_1=a(1,\,0,\,0)=-a_1+a_2+a_3,\quad A_2=a(0,\,1,\,0)=a_1-a_2+a_3,\\
A_3=a(0,\,0,\,1)=a_1+a_2-a_3.
\end{gathered}
\end{equation}
We let $\proj$ be the orthogonal projection of $\fcc$
onto the plane with normal given by $a_1-a_2=\frac a2(-1,\, 1,\, 0)$
(or the $\planev$ plane using Miller indices since we have that
$A_1-A_2=-2(a_1-a_2)$),
and let {$\qproj=I-\proj$} be the projection onto the line parallel to
$a_1-a_2=\frac a2(-1,\, 1,\, 0)$ (or the $\dirv$ direction using Miller indices~\cite{ashcroft}).
We observe that the projection of $\fcc$ onto
the plane normal to $a_1-a_2=\frac a2(-1,\, 1,\, 0)$ (or the $\planev$ plane using Miller indices)
is a triangular lattice $\tri := \proj \fcc$.  Each point in the triangular lattice
$\tri := \proj \fcc$ is thus the projection of a column of atoms with spacing $a/\sqrt 2$
parallel to
$a_1-a_2=\frac a2 (-1,\,1,\,0)$ (or the $\dirv$ direction
using Miller indices)
as depicted in Figure~\ref{fig: fcc projection figure}.

Now let $\pdomain$ be a finite subset of the triangular lattice $\tri$, and let
$\domain$ be the right cylinder over $\pdomain$ in the fcc lattice $\fcc$. We will
call a displacement $U:\domain \rightarrow \Real^3$ \emph{\periodic} if
\begin{equation*}
U(X) = U \bigg (X +  n(a_1-a_2)  \bigg) \quad\text{ for all } X \in \domain \text{ and all } n \in \Int.
\end{equation*}
We note that the scale of the periodicity is $|a_1-a_2|=a/\sqrt 2.$
Any periodic displacement $U$ reduces to
a displacement $\dsp: \pdomain \rightarrow \Real^3$ defined by
\begin{equation*}
\dsp(x) := U (X) \quad\text{ for any } X\in\domain \text{ with } \proj X = x.
\end{equation*}
We let $\Ub$ denote the space of displacements of $\pdomain$.

We will now derive an \emph{energy per period} $\energya : \Ub
\rightarrow \Real$ and a corresponding force field on $\pdomain.$ Let
$\p$ be a pair potential. For example, we use the Morse potential in
our numerical experiments below.  We will adopt a fourth-nearest
neighbor pair-interaction
model.  For this choice of the interaction model the quasi-nonlocal
coupling method~\cite{Shimokawa:2004} can be applied. {We note that,
throughout, we understand $n$-th nearest neighbours in terms of
Euclidean distance in the reference configuration in $\R^3$ (i.e., in the 3D FCC lattice)~\cite{Shimokawa:2004}.}

{Let $\N$ denote the $\fcc$-{\em interaction range}, by which
  we mean the set of all vectors pointing from an atom in $\fcc$ to
  one of its first, second, third, or fourth nearest-neighbors.
 Representative first, second, third, or fourth nearest-neighbor vectors are
  given in terms of the primitive lattice vectors~\eqref{prim} and cubic lattice vectors~\eqref{cubic} by
  \begin{equation}\label{firstsecond}
a_1= \left(0,\,\frac a2,\, \frac a2\right),\quad
A_1=\left(a,\,0,\,0\right),\quad
a_2+a_3= \left(a,\,\frac a2,\, \frac a2\right),\quad
2a_1=\left(0,\,a,\, a\right),
\end{equation}
with the full set of first, second, third, or fourth nearest-neighbor vectors
given by symmetry.
  }

For $X\in\domain,$ we let
\[
\N_X := \{B \in \N: X + B \in \domain \}.
\]
We will often call directions $B \in \N_X$ {\em bonds}.
We note that $\N\ne \N_X$ when $X$ lies in a shallow
layer along the boundary of $\domain.$
We also note that we have $\N_{X_1}=\N_{X_2}$ if
$PX_1=PX_2=x,$ hence we can unambiguously use the notation
$\N_x:=\N_X$ where $PX=x$ for $x\in\pdomain.$

 For any $B\in\N$ and any displacement $U(X),$ we
define the difference operator $\d_B$
by
\begin{equation*}
\d_B U(X) := U(X+B) - U(X),
\end{equation*}
whenever $X,\ X+ B \in \domain$.  We define the set of projections $b=PB$
of the bonds in $\N$ and $\N_X$ onto the plane $\planev$ by
\[
\M:=P\N\qquad\text{and}\qquad \M_x:=P\N_X.
\]
For periodic displacements $U(X),$ we define
\begin{equation*}
\d_b u(x):=u(x+b) - u(x) = U(X+ B) - U(X)=\d_B U(X)
\end{equation*}
where $u(x)=U(X)$ for $x=PX$ and $b=PB.$ 
It will also be convenient to use the notation
\[
\d_B u(x):=\d_b u(x)=\d_B U(X) \quad\text{for }PB=b.
\]

The \emph{energy per period} $\energya : \Ub \rightarrow \Real$ is
given by
\begin{equation*}
\energya(u) := \half \sum_{x \in \pdomain} \sum_{B \in \N_x}  \p(|B + \d_B u(x)|)
\end{equation*}
where we note that the difference in the deformed positions of atoms at
$X+B$ and $X$ is $B + \d_B u(x).$
The corresponding force at $x\in\pdomain$ is then given by
\begin{equation}\label{force}
\F^a(u)(x) := -\pder{}{u(x)}\energya(u)(x)=-\pder{}{u(x)}\sum_{B \in \N_x} \p(|B + \d_B u(x)|),
\end{equation}
where $\pder{}{u(x)}$ above denotes the
partial derivative with respect to the displacement $u(x)$ of the atom with projected reference
position $x\in\pdomain.$
We define the first variation (or, gradient) of $\energya(u)$ by $\del
\energya(u),$ 
which is given for each $x\in\pdomain$ by
\begin{equation}\label{first}
\del \energya(u)(x):=\pder{}{u(x)}\energya(u)(x).
\end{equation}
We can then express \eqref{force} in operator notation by
\[
\F^a(u)=-\del \energya(u).
\]

We will now formulate a boundary value problem for the energy
$\energya$. We let $\Gamma \subset \pdomain$ denote the part of
the ``boundary'' where the displacements are constrained. In
the numerical experiments described below,\ $\Gamma$ will be a
shallow layer of atoms along all or some of the boundary of
$\pdomain$.
Let the space of admissible displacements that are equal to
$\dsp_0\in \Ub$ on the boundary $\Gamma$ be denoted by
\begin{equation*}\label{admis}
\admis := \{\dsp \in \Ub: \dsp(x) = \dsp_0(x) \text{ for all } x \in \Gamma \}.
\end{equation*}
We will study the minimization problem
(or more precisely, the local minimization problem):
\begin{equation}\label{eqn: minimisation problem}
\text{Find } \dsp \in \newargmin_{v \in \admis} \energya(v).
\end{equation}
The Euler--Lagrange equation corresponding to problem~\eqref{eqn: minimisation problem}
is
\begin{alignat*}{3}
 -\del \energya(\dsp)(x) &
 = 0 &&\quad\text{ for all } x \in \pdomain \setminus \Gamma,\\
                  \dsp(x) &= \dsp_0(x) &&\quad\text{ for all } x \in \Gamma.
\end{alignat*}

In our numerical experiments, we consider only initial
displacements $u_0$ satisfying $u_0: \pdomain \rightarrow
\planev,$ where we recall that $\planev$ is the plane with
normal given by $a_1-a_2=\frac a2(-1,\, 1,\, 0)$. We will call
any displacement $u: \pdomain \rightarrow \planev$ a
\emph{plane displacement}, and we will let $\Ut\subset \Ub$
denote the space of plane displacements. For plane
displacements, $\energya$ can be interpreted as the energy of
the two dimensional crystal $\pdomain$ computed using a pair
potential which depends on the bond direction.

The restriction of the energy per period $\benergya(u):\Ub \rightarrow \Real$ to plane displacements $\Ut$
can then be given by
\begin{equation} \label{eqn: fully 2d energy}
\benergya(u) =\half \sum_{x \in \pdomain} \sum_{b \in \M_x}  \p_b(|b + \d_b u(x)|),\quad u\in\Ut,
\end{equation}
where $\p_b: (0,\infty) \rightarrow \Real$ is given by
\begin{equation*}
\p_b(r) := \sum_{\{B \in \N_x :\, \proj B = b \}} \p\left(\left(r^2 + |\qproj B|^2\right)^{1/2}\right)
\end{equation*}
{where we recall that $\qproj=I-P$ is the projection onto the line parallel to
$a_1-a_2=\frac a2(-1,\, 1,\, 0).$}
{Roughly speaking, $\p_b(|b + \d_b u(x)|)$ is the energy of the interaction of atom $x$ with
all of its neighbors in the column over $x+b.$}
For $x \in \pdomain$, we will call the inner sum
\begin{equation*}
\benergya_x(u) := \half\sum_{b \in \M_x}  \p_b(|b + \d_b u(x)|)
\end{equation*}
the \emph{atomistic energy} at $x$. 

It will be convenient to establish coordinates adapted to the
triangular lattice $\tri$, which lies in the plane defined by the
normal $a_1-a_2=\frac a2(-1,\, 1,\, 0)$ (or the $\planev$
plane).  This plane is spanned by the orthogonal vectors of
length $a/\sqrt 2$ in the $\hkl[1 1 0]$ direction and of length
$a$ in the $\hkl[0 0 1]$ direction given by
\begin{equation*}
\begin{split}
  V_1 :&= \smfrac{1}{2}(A_1+A_2)=a_3=\smfrac{a}{2}(1,\,1,\,0),  \\
  V_2 :&= A_3=a_1+a_2-a_3=a(0,\,0,\,1).
\end{split}
\end{equation*}
Throughout the remainder of this paper, we will denote the
coordinates of $(x_1,x_2)$ in the $\{\,V_1,\,V_2\,\}$ basis by
\begin{equation*} 
\<\nu_1,\nu_2\> := \nu_1 V_1 + \nu_2 V_2.
\end{equation*}
The triangular lattice $\tri$ is then generated by the basis
vectors $\<1,\,0\>$ and $\<1/2,\,1/2\>.$

It can be checked that there
are four distinct symmetry-related projected interaction potentials, $\p_b(r),$
corresponding to interactions with bonds $B$ such that
\begin{equation*}
\begin{alignedat}{2}
b_1&=\<1,\,0\>&&\quad\text{$PB=b_1$ and $QB=n(a_1-a_2)\ $ for $\ n=-1,\,0,\,1$},\\
b_2&=\<2,\,0\>&&\quad\text{$PB=b_2$ and $QB=0$},\\
b_4&=\<1/2,\,1/2\>&&\quad\text{$PB=b_4$ and $QB=n(a_1-a_2)\ $ for $\ n=-3/2,\,-1/2,\,1/2,\,3/2$},\\
b_5&=\<3/2,\,1/2\>&&\quad\text{$PB=b_5$ and $QB=n(a_1-a_2)\ $ for  $\ n=-1/2,\,1/2$.}
\end{alignedat}
\end{equation*}
where the planar bonds $b_i$ are displayed in Figure~\ref{fig: bonds figure}.
{The remaining interaction potentials, $\p_b(r)$ for
$b\in\M,$ can be obtained by symmetry from
${\p}_{b_1}(r)={\p}_{b_6}(r)={\p}_{b_7}(r),$
${\p}_{b_2}(r),$ ${\p}_{b_4}(r),$ and ${\p}_{b_5}(r).$
We note that the lattice constant $a$ is contained in the definition
given in the previous paragraph for the coordinates $\<\nu_1,\nu_2\>.$
}

Now observe that for $u \in \Ut$ we have $\del \benergya (u) \in \Ut$
where the gradient $\del \benergya (u)$ is defined by~\eqref{first}.
Therefore, if a gradient-based minimization algorithm (such as the steepest descent method or the nonlinear
conjugate gradient algorithm discussed below) is started from a plane displacement,
it will terminate at a plane displacement. Thus, since we will only consider initial
configurations which are plane displacements, we will be able to use the restriction
of $\benergya$ to $\Ut\cap\admis$ in our numerical experiments. This is important
since $\benergya$ restricted to $\Ut\cap\admis$ is the energy
of a two dimensional crystal, and in that case the patch-test consistent atomistic to continuum coupling method
developed in \cite{shapeev} can be used.

\subsection{Energy-based quasicontinuum approximations}
\begin{figure}
\centering
\includegraphics[width = 10cm]{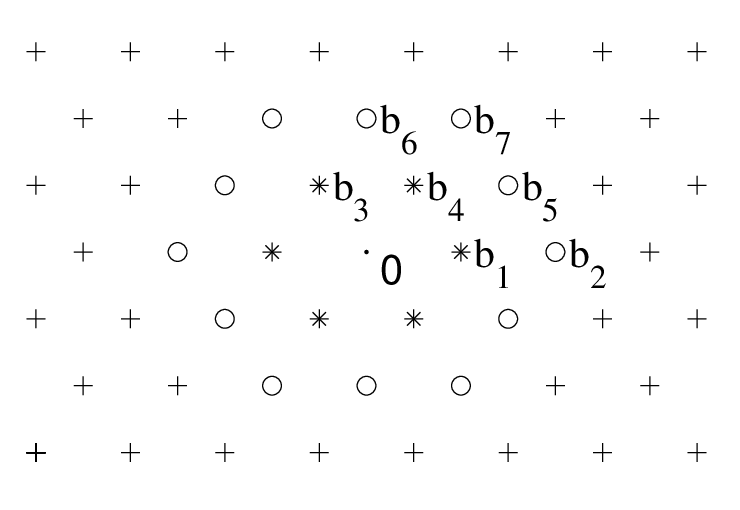}
 \caption{In the above display of the $\planev$ plane, the origin is marked with a dot, the nearest neighbors of the origin
	  are marked with asterisks, and the other neighbors are marked with circles.
	  The bond vector $b_i$ is the vector pointing from the origin to the atom marked $b_i$.}
\label{fig: bonds figure}
\end{figure}

We will construct a local approximation to the energy $\benergya$. Our
methods follow \cite{Shimokawa:2004}. First, we define an extrapolated
difference operator $\dextra$. The difference operator $\dextra$ will
approximate the vector pointing from an atom to one of its neighbors
using only the vectors pointing from an atom to its nearest
neighbors. The \emph{nearest neighbors} are the atoms depicted in
Figure~\ref{fig: bonds figure}. {Recall that $\M := \proj \N$; we thus have $b_i \in
\M$ for the vectors $b_i$, $i = 1, \dots, 7$, depicted in
Figure~\ref{fig: bonds figure}, and moreover, $\M$ is obtained by
applying the lattice symmetries of $\tri$ to $b_1, \dots,
b_7$. For these vectors, we define the operators
\begin{align*}
\dextra_{b_i} :&= \d_{b_i} \text{ for } i \in \{1, 3, 4\}, \\
\dextra_{b_2} :&= 2\d_{b_1}, \\
\dextra_{b_5} :&= \d_{b_1} + \d_{b_4},\\
\dextra_{b_6} :&= \d_{b_3} + \d_{b_5}, \\
\dextra_{b_7} :&= 2\d_{b_4},
\end{align*}
and extend the definition by symmetries of $\tri$ so that $\dextra_b$
is defined for all $b \in \M$.}

Using the difference operator $\dextra$,
we define the \emph{Cauchy--Born energy} by
\begin{equation*} 
\E^{cb}(u) := \half\sum_{x \in \pdomain} \sum_{b \in \M_x}  \p_b(|b + \dextra_b u(x)|).
\end{equation*}
{We will call the inner sum
\begin{equation*}
\E^{cb}_x(u) :=  \half\sum_{b \in \M_x}  \p_b(|b + \dextra_b u(x)|)
\end{equation*}
the \emph{continuum energy} at x,
so
\[
\E^{cb}(u)= \sum_{x \in \pdomain} \E^{cb}_x(u).
\]
The following remark gives a more
detailed motivation for these definitions.
}

{
\begin{remark}[The Cauchy--Born Approximation]
  We call $\E^{cb}(u) $ the Cauchy--Born energy since, if the
  displacement $u$ is interpreted as a piecewise linear spline with
  respect to the canonical triangulation of $\tri$, it can be
  rewritten as an integral over a stored energy density.

  To make this precise, we observe that the continuum energy
  at $x\in\pdomain$ such that $x+b\in\pdomain$ for all $b \in \M$
   is given for uniform displacements $u^F(x):=Fx-x$ where $F\in \R^{2 \times 2}$ by
\[
\E^{cb}_x\big(u^F\big)=\half \sum_{b \in \M}\p_b(|Fb|).
\]
We thus
define the {\em Cauchy--Born stored
    energy density,} $W_{cb} : \R^{2 \times 2} \to \R \cup
  \{+\infty\}$, by
\[
W_{cb}(F)=\frac{\half\sum_{b \in \M}\p_b(|Fb|)}{v},
\]
  where $v$ is the area associated with each atom (in the 2D
  lattice). Moreover, for $u \in \Us$, let $\bar{u}$ denote the piecewise
  affine interpolant of $u$ with respect to the canonical
  triangulation of $\tri$, let $\nabla\bar{u}$ denote the resulting
  displacement gradient, and let $\bar{\omega}$ denote the union of
  those elements.

For periodic boundary conditions, or by modifying the functional
$\E^{cb}(u)$ at the boundary $\Gamma$ of the domain $\pdomain,$ one can use
Shapeev's
  bond density lemma \cite{shapeev} to prove that
  \begin{displaymath}
    \E^{cb}(u) = \int_{\bar{\omega}} W_{cb}(\nabla \bar{u}) \, {\rm dx}.
  \end{displaymath}
  We note that this modification of the functional
$\E^{cb}(u)$ at the boundary does not affect the minimization of quasicontinuum energies
whose set of admissible
displacements, $\admis,$ are constrained in the entire boundary
of $\pdomain,$ such as for the Lomer Dislocation Dipole problem that we
study.
  More detailed discussions and analyses of the Cauchy--Born
  approximation can be found in
  \cite{FrTh:2002,BLBL:arma2002,E:2007a}.
\end{remark}
}

Now let $\pdomain := \A \cup \C$ be a partition of $\pdomain$ into an
\emph{atomistic region} $\A$ and a \emph{continuum region} $\C$.  We
define the \emph{energy based quasicontinuum } (QCE) energy by
\begin{equation} \label{eqn: qce energy}
\E^{qce}(u) := \sum_{x \in \atomistic} \benergya_x(u)
                    + \sum_{x \in \continuum} \E^{cb}_x(u).
\end{equation}

Following \cite{Shimokawa:2004}, we also define the \emph{quasi-nonlocal } (QNL) energy by
\begin{equation}\label{eqn: qnl energy}
\E^{qnl}(u) := \half\sum_{x \in \pdomain} \sum_{b \in \M_x}
\bondind(x) \p_b(|b + \dextra_b u(x)|) + \left(1 - \bondind(x)\right) \p_b(|b + \d_b u(x)|),
\end{equation}
where
\begin{equation*}
\bondind(x) :=
\cases{
1 &\quad\text{ if } x\in\continuum \text{ and } x+b \in \continuum, \\
0 &\quad\text{ otherwise.}
}
\end{equation*}

{ The reason for the introduction of the QNL method was that
  the QCE method did not pass the {\em patch test}
  \cite{Shenoy:1999a,Dobson:2008b}, that is, the coupling mechanism
  defined in \eqref{eqn: qce energy} results in non-zero forces (the
  ``ghost forces'') at the atomistic/continuum interface under
  homogeneous displacements (or, deformations). By contrast, it was
  shown in \cite{Shimokawa:2004} that the QNL energy \eqref{eqn: qnl
    energy} does pass the patch test \cite{Shimokawa:2004,
    Dobson:2008b,doblusort:qce.stab}, as long as two atoms interact
  only when they share a common nearest neighbor. The authors of
  \cite{Shimokawa:2004} consider only the case where the pair
  potential does not depend on the bond direction. Nonetheless, the
  argument which they give may be applied to show that the QNL
  energy~\eqref{eqn: qnl energy} passes the patch test also in the
  present case.

  The importance of passing the patch test lies in the fact that the
  ``ghost forces'' can be understood as a consistency error, which
  results in an $O(1)$ relative error in the displacement gradient
  \cite{Dobson:2008b}. Moreover, it was shown in
  \cite{doblusort:qce.stab}, that they can result in a $O(1)$
  relative error in the prediction of critical loads at which lattice
  instabilities occur.
}

To achieve an efficient quasicontinuum method, the positions of atoms
in the continuum region are normally further constrained by piecewise
linear interpolation~\cite{Miller:2003a}.  The development of an
implementable and patch-test consistent quasicontinuum energy that
allows coarsening by piecewise linear interpolation in the continuum
region has been achieved so far only for two-dimensional problems with
pair potential interactions~\cite{shapeev}.  The development of
patch-test consistent quasicontinuum energies for many-body potentials
such as the popular embedded atom method (EAM) or for
three-dimensional problems is an open problem.

\subsection{The force based quasicontinuum approximation}

The force based quasicontinuum (QCF) approximation gives a
patch-test consistent approximation for atomistic-to-continuum
coupling even with coarsening in the continuum
region~\cite{curt03,dobsonluskin07}. Force-based multiphysics
coupling methods are generally popular because of their
algorithmic simplicity. However, force-based coupling methods
are known to give non-conservative force
fields~\cite{curt03,dobsonluskin07}, and our recent research
described at the end of this section has discovered additional
stability problems for QCF approximations and solution
methods~\cite{qcf.iterative}.

Instead of approximating the energy $\benergya$ by a local energy, we can approximate the forces $\F(u)$ directly. This
leads to the \emph{force based quasicontinuum } (QCF) approximation.
As above, let $\pdomain := \A \cup \C$ be a partition of $\pdomain$ into an atomistic region $\A$
and a continuum region $\C$. Define $\forceqcf: \Ut \rightarrow \Ut$
by
\begin{equation}  \label{eqn:qcf_force}
 \forceqcf(u)(x) :=
\cases{
- \del \benergya(u)(x) &\quad\text{ if } x \in \A, \\
- \del \E^{cb}(u)(x) &\quad\text{ if } x \in \C.}
\end{equation}
We consider the boundary value problem
\begin{equation}\label{eqn: qcf bvp}
\begin{alignedat}{3}
 \forceqcf(u)(x) &= 0 &&\quad\text{ for all } x\in \pdomain \setminus \Gamma, \\
            u(x) &= u_0(x) &&\quad\text{ for all } x \in \Gamma.
\end{alignedat}
\end{equation}
It was shown in
  \cite{dobs-qcf2,doblusort:qcf.stab,MakrOrtSul}, for one-dimensional
  model problems, that the solution of problem~\eqref{eqn: qcf bvp}
approximates the solution of problem~\eqref{eqn: minimisation
  problem}.  It is also true that the QCF method passes the patch test
\cite{dobsonluskin07}.

In our numerical experiments we will solve equation~\eqref{eqn: qcf bvp}
 using an iterative method defined in \cite{dobsonluskin07,dobsonluskin08}.
Let $u_0 \in \Ut$ be given, and let $\E^{qc}$ be either $\E^{qnl}$ or $\E^{qce}$.
 Then we let $u_{n+1} \in \admis \cap \Ut$ be the solution of the problem:
\begin{equation}\label{eqn: gfc iteration}
 u_{n+1} \in \newargmin_{v \in \admis \cap \Ut} \left \{
 \E^{qc}(v) - \sum_{x \in \pdomain} \left(\forceqcf(u_n)(x) + \del \E^{qc}(u_n)(x)\right) \cdot v(x) \right \}.
\end{equation}
The method \eqref{eqn: gfc iteration} with $\E^{qc}=\E^{qce}$ is the
popular {\em ghost force correction} method proposed in
\cite{Miller:2003a}. We have proposed the method $\E^{qc}=\E^{qnl}$ as
a {\em ghost force correction} method with improved stability
properties~\cite{qcf.iterative}. {For $\E^{qc}=\E^{qce},$ we call the
iteration \eqref{eqn: gfc iteration} the QCF-QCE iteration; for
$\E^{qc}=\E^{qnl},$ we call it the QCF-QNL iteration.}

The Euler--Lagrange equation for problem~\eqref{eqn: gfc iteration} is
 \begin{alignat*}{3}
   \del \E^{qc}(u_{n+1})(x) &= \forceqcf(u_n)(x) + \del \E^{qc}(u_n)(x) &&\quad\text{ for all } x \in \pdomain \setminus \Gamma, \\
  u_{n+1} (x) &= u_0(x) &&\quad\text{ for all } x \in \Gamma.
 \end{alignat*}
It is easy to see that, if the sequence of iterates converges,
then its limit must be a solution of problem~\eqref{eqn: qcf
bvp}. Moreover, under certain technical conditions related to
the stability of the preconditioner $\E^{qc},$ it can be shown
that the sequence indeed converges
\cite{dobsonluskin07,dobsonluskin08,qcf.iterative}.

A QCF solution method that can theoretically give inaccurate
results is the use of a nonlinear conjugate gradient solver for
the force-based quasicontinuum
equations~\cite{Miller:2008,cadd2,curt03}. We have proven that
the linearization of the force-based equilibrium equations is
not positive definite~\cite{dobs-qcf2}, which implies that the
conjugate gradient solution of this problem is
unstable~\cite{qcf.iterative}.  We have discovered from
informal discussions with computational physicists that their
conjugate gradient iterative solution of force-based
multiphysics coupling methods sometimes oscillates rather than
converges, a phenomenon partially explained by our theoretical
analysis of this instability. We are developing benchmark tests
to further study the reliability of the conjugate gradient
solution of the force-based quasicontinuum equations.

We have shown theoretically and computationally for one-dimensional
model problems that the GMRES method is a reliable and efficient
solver for the force-based quasicontinuum
equations~\cite{qcf.iterative}. We are thus also preparing
multi-dimensional benchmark tests to evaluate the reliability and
efficiency of the GMRES solution of the force-based quasicontinuum
equations.

\section{Numerical Experiments}
\label{sec:numerical_experiments}

\subsection{The Lomer Dislocation Dipole: Analytical Model}

\begin{figure}\centering
  \includegraphics[width = 11.5 cm]{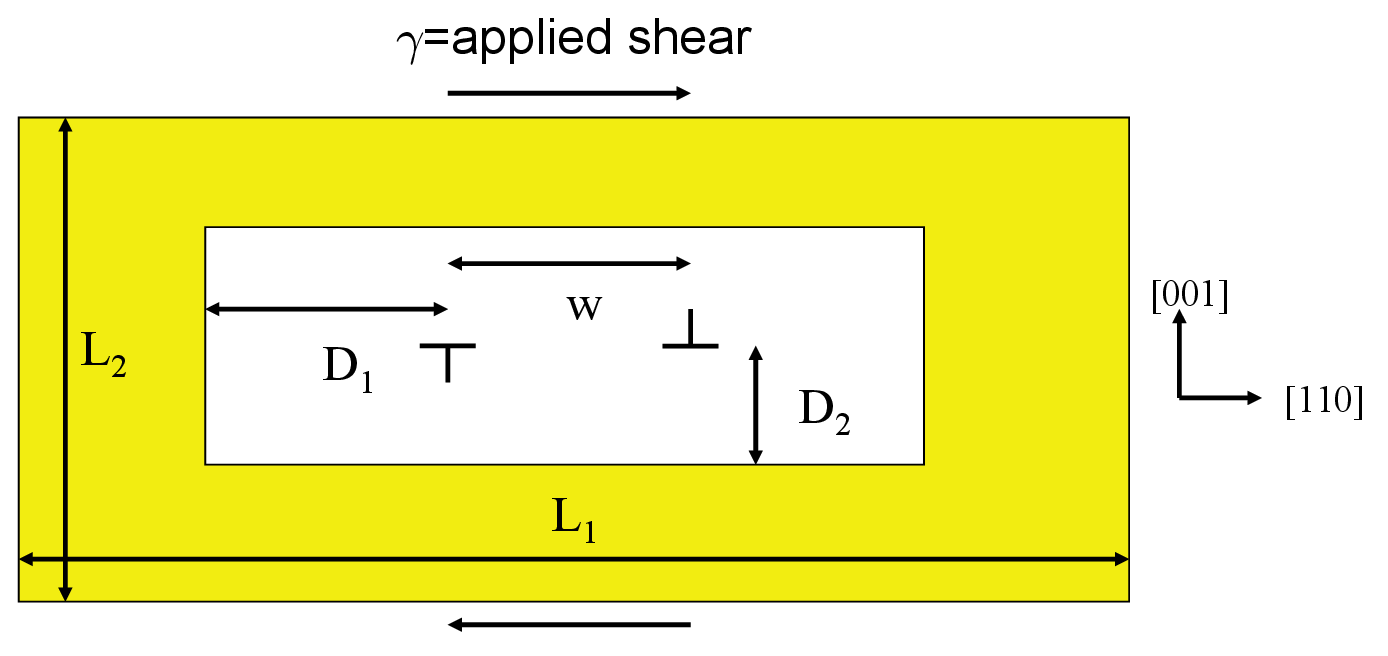}
  \caption{The atomistic and continuum (yellow) domains for the Lomer dislocation dipole.
  The size of the total computational domain is $L_1=75\times a/\sqrt 2$ and
  $L_2=60\times a$ where $a$ is the length of the sides of the cubic supercell.
  The atomistic domain is given by $D_1= k \times a/\sqrt 2$
  and $D_2= k \times a$ for $k=3,\,4,\,5.$ The separation of the dipole used in
  the initial guess $u_{elas}$ is
  $w=10\times a/\sqrt 2.$ }
  \label{dipole}
\end{figure}

Following the numerical experiments presented in~\cite{Miller:2008}, we consider
a dipole of Lomer dislocations, as depicted in Figure~\ref{dipole}. The Lomer dislocation
has Burgers vector $b=\hkl[1 1 0]$ in the $\hkl(0 0 1)$ plane~\cite{hirthlothe}.
The dipole should be a stable equilibrium
if the distance between the cores is large enough so that the
Peach-Koehler elastic attraction of the dislocations is weaker than the critical
force needed to overcome the Peierls energy barrier~\cite{Miller:2008}.

We will study the stability of the dipole
when a shear is applied to
the crystal. By a shear, we mean a homogeneous deformation of the
lattice $\tri$ that takes the form
\begin{equation*}
x \mapsto \sigma(\gamma)x\quad\text{where}\quad\sigma(\gamma):=\left (
    \begin{array}{cc}
      1 & \gamma \\
      0 & 1
    \end{array} \right ).
\end{equation*}
We call $\gamma$ the \emph{shear strain}.
We say that the shear is \emph{positive} if $\gamma > 0$ and
\emph{negative} if $\gamma < 0$. More generally, we will apply
a shear to the deformation $x+u(x)$ given by the displacement
$u(x)$ to obtain the sheared deformation
\begin{equation*} 
x+u(x) \mapsto \sigma(\gamma)\left(x+u(x)\right).
\end{equation*}

A positive shear applied to the top and bottom boundaries of
the crystal changes the energy landscape to favor the movement
of the dislocations apart. As $\gamma$ increases, the shear
force overwhelms both the Peierls force and the Peach-Koehler
elastic attraction. Thus, the original equilibrium
configuration of the dipole becomes unstable, and the
dislocations move away from each other. We call the value of
$\gamma$ at which this instability occurs the \emph{critical
strain}. In our numerical experiments, we will simulate the
process of slowly shearing the crystal until the critical
strain is reached.

To understand the movement of the Lomer dipole theoretically,
we model the local minima of the atomistic energy for
displacements constrained on the boundary $\Gamma$ to
be $\sigma(\gamma)x -x$ and constrained to be in the energy well of a Lomer dipole
with separation $w:$
\begin{align*}
\inf_{v\in\widetilde{{\rm Adm}}\left(\sigma(\gamma)x -x,\,w\right)}\E^a(v)\approx&\tilde\E^a(w,\gamma)= \E^a_{\text {misfit}}(w)+
\E^a_{\text {dipole attraction}}(w)\\
& +\E^a_{\text {applied shear}}(w, \gamma)
+\E^a_{\text {boundary effect}}(w, \gamma)
\end{align*}
where $\widetilde{{\rm Adm}}\left(\sigma(\gamma)x -x,\,w\right)$ is the set
of admissible displacements roughly described above. Here the classical
Peierls-Nabarro misfit energy~\cite{hirthlothe}, the classical Peach-Koehler interaction
energy~\cite{hirthlothe}, and our modeling of the effect of the applied shear are given for $w\ge
0$ and $\gamma\ge 0$ by
\begin{equation}\label{energies}
\begin{split}
&\E^a_{\text {misfit}}(w)=\frac{\mu b}{2\pi(1-\nu)}\cos\frac{2\pi w}b\exp\left[{-\frac{\pi d}{(1-\nu)b}}\right],\\
&\E^a_{\text {dipole attraction}}(w)=\frac{\mu b^2}{2\pi}\ln \frac w{L},\\
&\E^a_{\text {applied shear}}(w, \gamma)=\beta_1 w+\beta_2\gamma-\frac {\beta_3\gamma\, w}{L},
\end{split}
\end{equation}
where $\mu$ is a characteristic shear modulus, $\nu$ is a
characteristic Poisson's ratio, $L$ is a characteristic length,
$b=a/\sqrt 2$ is the magnitude of the Burger's vector,
 $d=a/\sqrt 2$ is the interatomic spacing between the
$\planev$ planes, and $\beta_1,\ \beta_2$ and $\beta_3$ are
constants with $\beta_3>0.$ Although our model for the misfit energy
$\E^a_{\text {misfit}}(w)$ and the dipole attraction
$\E^a_{\text {dipole attraction}}(w)$ are well-known
approximations of the atomistic energy $\E^a(v)$
\cite{hirthlothe}, our bilinear model for the shear energy $\E^a_{\text
{shear}}(w, \gamma)$ is not derived from the atomistic
energy and our justification is based only on the qualitative
behavior it predicts for the stability of the Lomer dipole.

\begin{figure}
  \includegraphics[width = 11.5 cm]{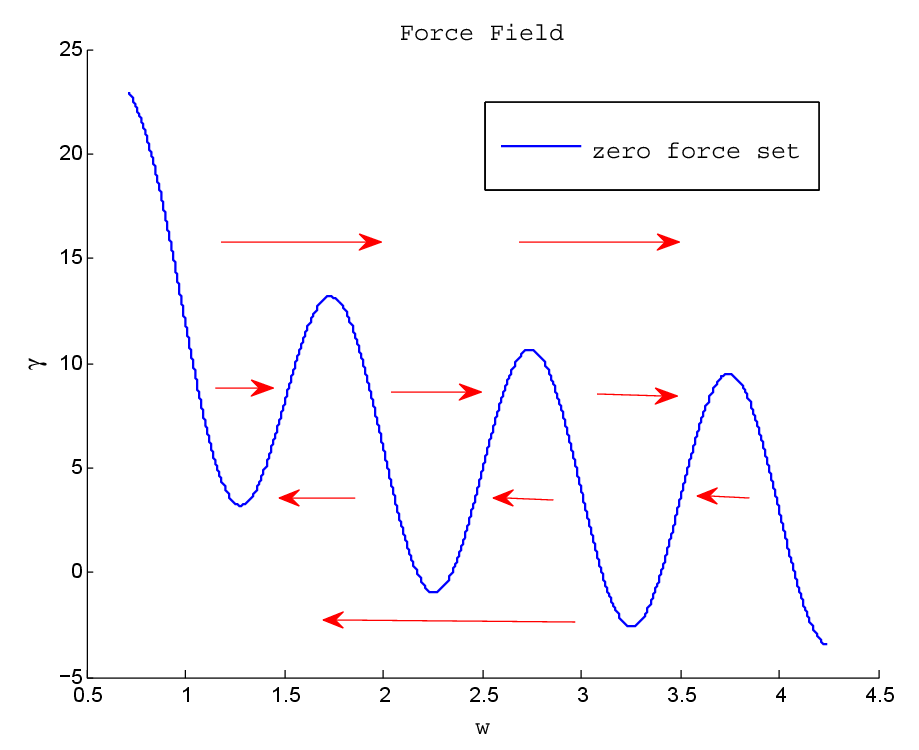}
\caption{The force field ~\eqref{conjugate} for the Lomer dipole model~\eqref{scale} with $\beta=12.$}
\label{force.field}
\end{figure}

\begin{figure}
  \includegraphics[trim=4.5cm 0cm 3.3cm 0cm, clip, width = 11.5 cm]{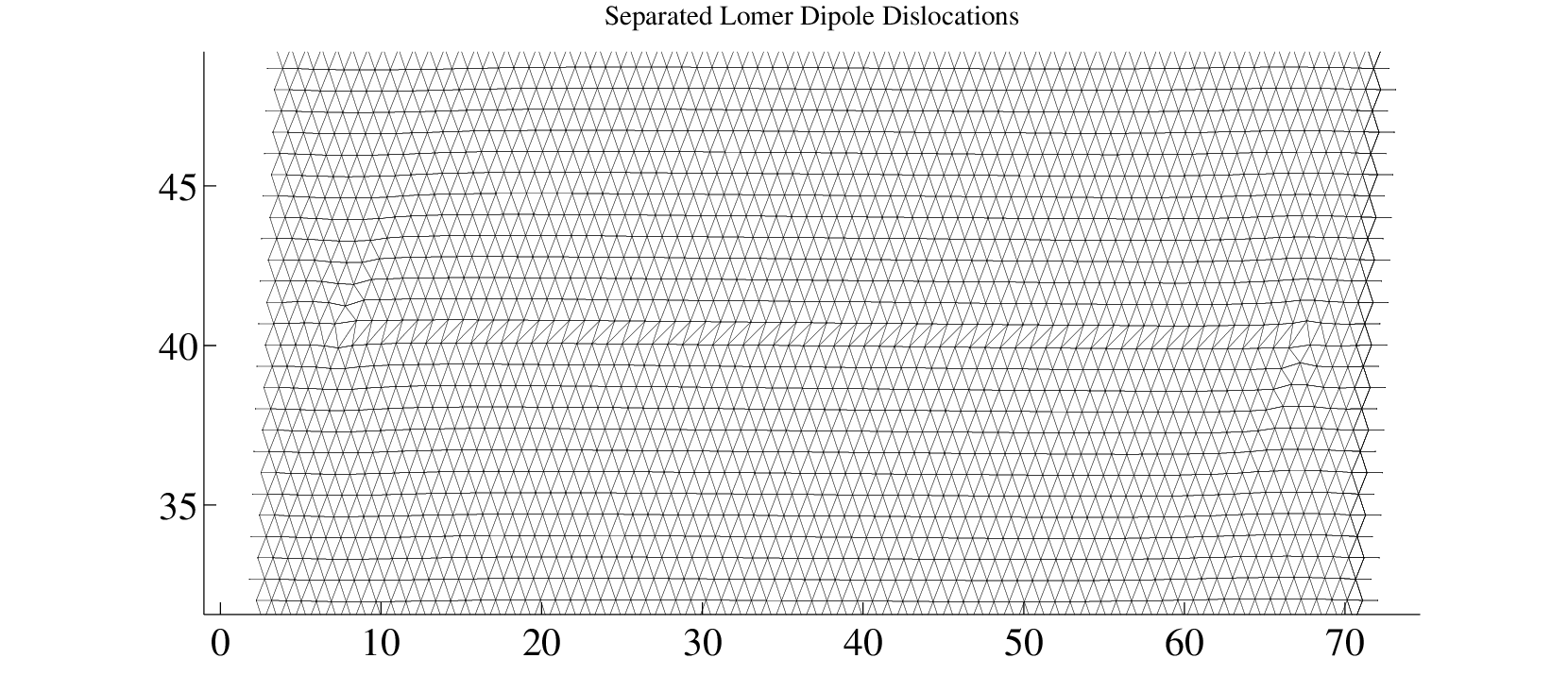}
\caption{Separation of the Lomer dislocations to the boundary after continuation past a critical strain
to $\gamma=0.070.$ Cartensian coordinates are
used to label the axes.  The top and bottom of the deformation of the reference domain $\pdomain$ are not displayed, but the
entire width of the deformed reference domain is displayed. The triangular
mesh was constructed from the reference lattice.}
\label{separate}
\end{figure}

After scaling the atomistic energy $\E^a(v),$ the dipole
separation $w,$ and the shear strain $\gamma,$ and after neglecting
the boundary effect $\E^a_{\text {boundary effect}}(w, \gamma)=0,$ we can obtain
from \eqref{energies} the model
\begin{equation}\label{scale}
\tilde\E^a(w,\gamma)=\cos w +\beta \ln w -\gamma w
\end{equation}
for some $\beta>0.$  We chose $\beta=12$ in our numerical experiments.
 We then find that the force conjugate to the dipole separation $w$ is
given by
\begin{equation}\label{conjugate}
-\pder{\tilde\E^a(w,\gamma)}{w} =\sin w - \frac {\beta} w +\gamma .
\end{equation}
We can see from the  force field displayed in
Figure~\ref{force.field} that the dipole separation $w$ becomes
unstable at a critical shear strain $\gamma$ and will tend to
infinity under gradient flow dynamics for the force field
~\eqref{conjugate} {(for example, consider the stability of the equilibrium
solution branch starting at
$w=1.5$ as $\gamma$
is increased)}.  The Lomer dislocations similarly separated
to the boundary at a critical shear strain $\gamma$ in our
numerical experiments (see Figure~\ref{separate}).

\subsection{The Lomer Dislocation Dipole: Numerical Experiments}

In out numerical experiments we use the Morse interatomic potential
\begin{equation*}
\phi(r) := \left(1 - \exp(\alpha(r - 1))\right)^2\quad\text{for }r>0,
\end{equation*}
with
\begin{equation*}
\alpha := 4.4.
\end{equation*}
We remark that we did not choose $\alpha = 4.4$ in order to model some
specific material. Rather, we chose this value of $\alpha$ as small as
possible while retaining a stable dipole.  For smaller values of
$\alpha$ we were unable to compute a stable dipole, whereas, if
$\alpha$ is chosen too large then the atomistic, Cauchy--Born, and
quasicontinuum models are almost identical since the first neighbor
interactions, which would become dominant, are treated identically in
all models.


We chose the lattice spacing $a$ so that
the FCC lattice $\mathcal{L}$ is a ground state for the atomistic energy
defined by the Morse potential and using a cut-off radius that includes
first, second, third, and fourth-nearest neighbor
 interactions (determined from the reference positions). 
Specifically, we let $a$ be the minimizer of
\begin{equation*}
\psi(r) := 12\phi \left ( \sqrt{\frac{1}{2}}r \right ) + 6\phi \left ( r \right )
    + 24\phi \left (\sqrt{\frac{3}{2}}r \right ) + 12\phi \left (\sqrt{2}r \right ),
\end{equation*}
which is the energy per atom in the undeformed lattice with lattice spacing $r$
since each atom in the $\fcc$ lattice has $12$ nearest neighbors, $6$ next-nearest
neighbors, $24$ third-nearest neighbors, and $12$ fourth-nearest neighbors
(see \eqref{firstsecond}).
The minimum is attained when $a = 1.3338$.

We let $\pdomain$ be the rectangular domain consisting of $150$
columns of atoms parallel to $V_2$ with each column containing
$60$ atoms defined in the $\<\nu_1,\,\nu_2\>$ coordinates by
\begin{equation*}
\begin{split}
\pdomain&=\left\{\<k,\,l\>: 0\le k\le 74,\ 0\le l\le 59\right\}\\
&\qquad\qquad\cup
\left\{\<k+1/2,\,l+1/2\>: 0\le k\le 74,\ 0\le l\le 59\right\}.
\end{split}
\end{equation*}
We take the boundary $\Gamma$ of $\pdomain$ to be a layer of
atoms four rows deep around the edges of $\pdomain,$ so that
\begin{equation*}
\begin{split}
\Gamma&=\left\{\<k,\,l\>: k=0,\,1,\,73,\,74,\ 0\le l\le 59\right\}\\
&\qquad\qquad\cup
\left\{\<k+1/2,\,l+1/2\>: k=0,\,1,\,73,\,74,\ 0\le l\le 59\right\}\\
&\qquad\qquad\cup
\left\{\<k+1/2,\,l+1/2\>: 0\le k\le 74,\ l=0,\,1,\,58,\,59\right\}\\
&\qquad\qquad\cup
\left\{\<k+1/2,\,l+1/2\>: 0\le k\le 74,\ l=0,\,1,\,58,\,59\right\}.
\end{split}
\end{equation*}
With this choice of $\Gamma$, for any atom $x \in \pdomain
\setminus \Gamma$, all fourth-neighbors of $x$ belong to
$\pdomain$. Thus, there are no boundary effects for the
reference lattice when $a = 1.3338.$  This choice of $\Gamma$
corresponds to taking Dirichlet boundary conditions on the
entire boundary. Throughout the rest of this section, we will
consider only displacements $u$ so that for some $\gamma>0$, we
have $u(x) = \sigma(\gamma)x-x$ for all $x\in\Gamma$. For such
a displacement, we will call the shear strain $\gamma$ the
\emph{shear on the boundary}.

\begin{figure}
\centering
 \subfigure[Atomistic region (3)]{\includegraphics[width = 8 cm]{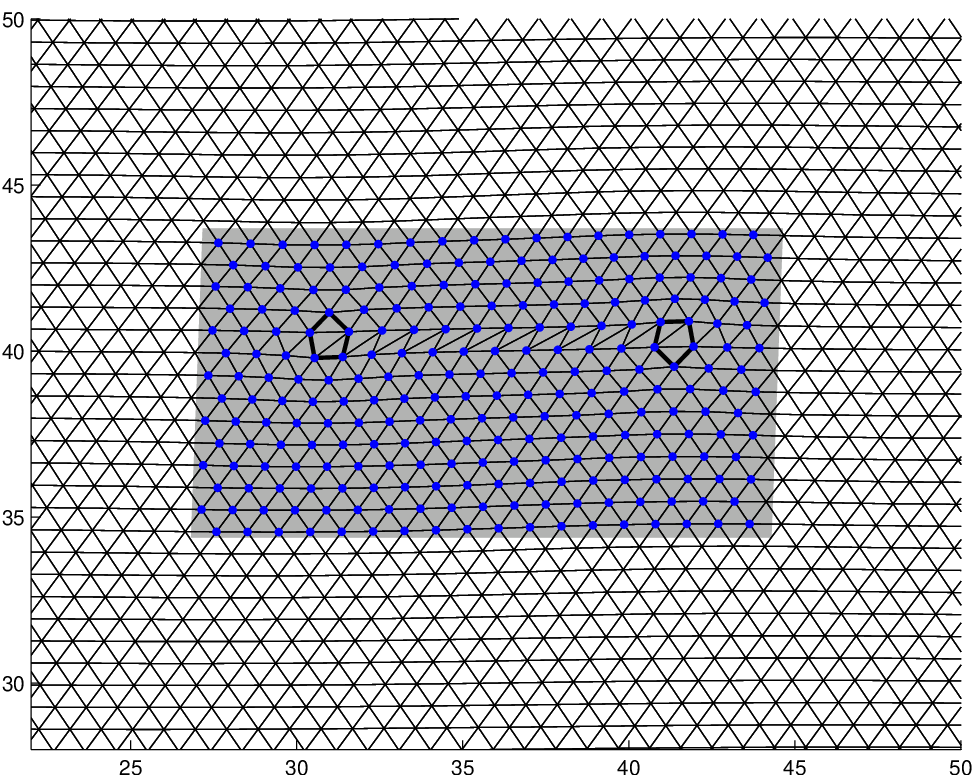}}
 \subfigure[Atomistic region (4)]{\includegraphics[width = 8 cm]{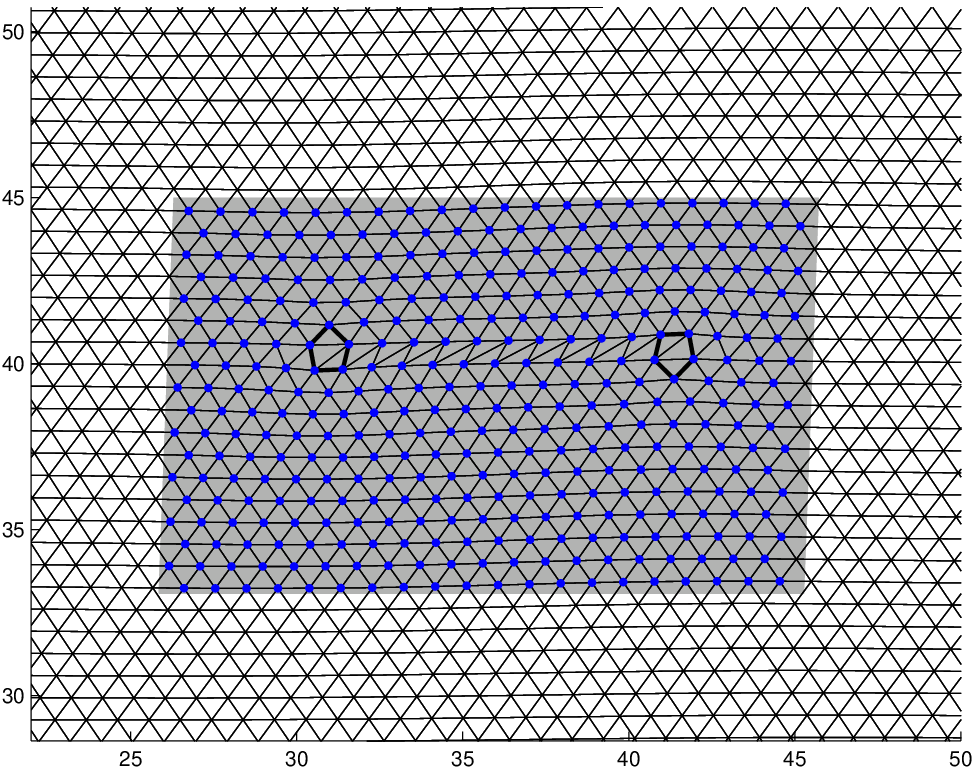}}
 \subfigure[Atomistic region (5)]{\includegraphics[width = 8 cm]{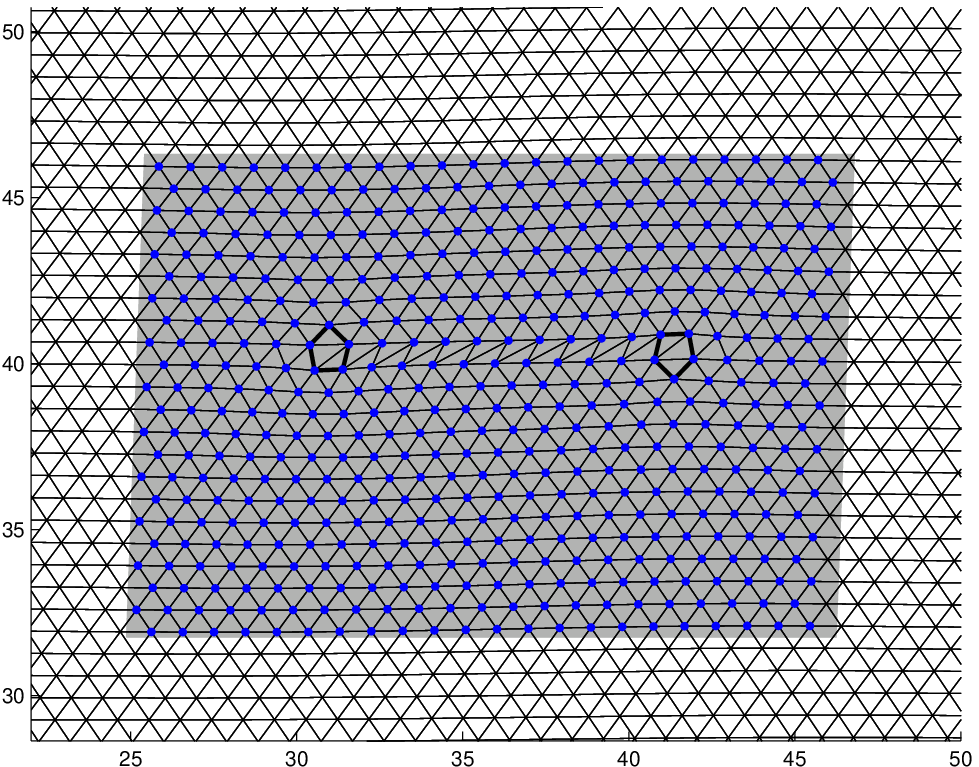}}
 \caption{The atomistic regions for the Lomer dislocation dipole.  Not
   all of the continuum region is displayed. The mesh is the Delaunay
   triangulation \cite{CompGeom} of the reference lattice (i.e., with
   the atoms as the nodes). The atomistic region is shaded grey. }
 \label{fig: dipole atomistic regions}
\end{figure}

We will now discuss how the stable Lomer dipoles displayed in
Figure~\ref{fig: dipole atomistic regions} are computed.  We first
construct an initial guess $u_{elas}$ for the displacement field
corresponding to a Lomer dipole by using the displacement fields of
the isotropic, linear elastic solution for the edge dislocations
\cite{hirthlothe} at the left and right pole of the dipole,
respectively. More precisely, the formula for these displacement
fields is
\begin{align*}
u_{1}^{*}=&\frac{b^{*}}{2\pi}
\Bigg[\tan^{-1}\frac{x_2-x_2^{*}}{x_1-x_1^{*}}+\frac{(x_1-x_1^*)(x_2-x_2*)}{2(1-\nu)\Big((x_1-x_1^*)^2+(x_2-x_2^*)^2\Big)}\Bigg],\label{EdgeDislocation}\\
u_{2}^{*}=&-\frac{b^*}{2\pi}\Bigg[\frac{1-2\nu}{4(1-\nu)}\ln\Big(x_1-x_1^*)^2+(x_2-x_2^*)^2\Big)\\
&\qquad\qquad+\frac{(x_1-x_1^*)^2-(x_2-x_2^*)^2}{4(1-\nu)\Big((x_1-x_1^*)^2+(x_2-x_2^*)^2\Big)}\Bigg],\nonumber
\end{align*}
where $*=L,R$ refers to the left or the right edge dislocation,
respectively, $b^*$ is the corresponding Burgers vector, and $x_i^*$,
$i = 1,2$, are the components of the position vectors~\cite{hirthlothe}. In our
numerical experiments, we set the initial positions of the Lomer
dislocations in the reference configuration to be (in the
$\<\nu_1,\,\nu_2\>$ coordinates)
\begin{equation*}
\left\<\nu_1^{L},\nu_2^{L}\right\>=\left\<32, 30+\frac{1}{6}\right\>\quad \text{and}
\quad \left\<\nu_1^{R},\nu_2^{R}\right\>=\left\<42, 30+\frac{1}{3}\right\>.
\end{equation*}
We note that initial dislocation positions are placed between
the center row of atoms at $\nu_2=30$ and the first row above
the center $\nu_2=30+1/2,$ with $\nu_2^L=30+1/6$ and
$\nu_2^R=30+1/3$ placed symmetrically about $\nu_2=30+1/4$ (see
Figure~\ref{fig: dipole atomistic regions}).

{ The displacements $u^{L}, u^{R}$, are estimates of the
  displacement fields for isolated edge dislocations. By superimposing
  them, we obtain an estimate of the elastic displacement field for
  the dipole (without shear):
\[
  u_{elas} = u^L + u^R.
\]
Finally, we apply a shear deformation to the deformation field $x + u_{elas}$,
\begin{equation*} 
 x+u^0_{elas}(x)= \sigma(\gamma) \left(x+u_{elas}(x)\right),
\end{equation*}
which yields an initial guess for the displacement field for the Lomer
dipole under applied shear:
\[
  u^{0}_{elas}(x)=\sigma(\gamma)\left(x+u_{elas}(x)\right)-x.
\]
}

We use the displacement $u_{elas}^0 = \sigma(\gamma_0) u_{elas}(x)-x$
as a starting guess for a preconditioned nonlinear conjugate gradient
method (P-nCG), which is described in Section \ref{sec:pcg}, to solve
the minimization problem~\eqref{eqn: minimisation problem}.

It is challenging to find a stable dipole. In fact, we were
unable to find a stable dipole without applying a positive
shear to the boundary. That is, we did not find a Lomer
dipole that was a local minimum of the atomistic
energy~\eqref{eqn: fully 2d energy} subject to the boundary
condition $u(x) = 0$ for all $x \in \Gamma$.  We were also
unable to find a stable dipole when the parameter $\alpha$ in
the Morse potential was too small. Essentially, we found that
as $\alpha$ decreased, the interval of shear strains for which
a given dipole was stable became smaller. After some
experimentation, we did find a stable Lomer dipole with $\alpha
= 4.4$ and with shear on the boundary $\gamma_0 = 0.0375$.

{Our atomistic and quasicontinuum models utilize a fourth-nearest
neighbor cutoff calculated from the positions in the reference
lattice, which is in fact the largest possible neighborhood for the
QNL method applied to an FCC lattice~\cite{Shimokawa:2004}.  Such a
cutoff is acceptable for atomistic configurations that are ``close''
to the reference lattice, however, for the highly deformed positions
near the dipole a cutoff in {\em deformed} coordinates ought to be
used. Even though the dominant nearest-neighbour interactions are
always included in our calculation, a more precise understanding of
the error committed for second- to fourth-nearest neighbours is
required.}

 Next we describe how we simulate the shearing of the dipole.
For simplicity, we will explain how this is done for
energy-based methods before we discuss the force-based method.
Let $u_0$ be the displacement field of the stable Lomer dipole
discussed above. Suppose that $\gamma_0$ is the shear on the
boundary for $u_0$.  Let $\gamma_k$ be an increasing sequence
of shear strains, and let $\del \gamma_{k+1} : = \gamma_{k+1} -
\gamma_k$.  $\mathcal{E}$ can be either the atomistic
energy~\eqref{eqn: fully 2d energy},\ the QCE energy
\eqref{eqn: qce energy}, or the QNL energy\ \eqref{eqn: qnl
  energy}.

We perform the following iteration starting with the displacement $u_0$.
Suppose we have a dipole with displacement field $u_k$
which solves the boundary value problem
\begin{equation*}
\begin{alignedat}{3}
\mathcal{- \del \mathcal{E}}(u_{k})(x) &= 0  &&\quad\mbox{ for all } x \in \pdomain \setminus \Gamma, \\
u_{k} (x) &= \sigma(\gamma_{k})x-x &&\quad\mbox{ for all } x
\in \Gamma. \nonumber
\end{alignedat}
\end{equation*}
We let
\begin{equation*}
u^{0}_{k+1}(x) = \sigma(\del \gamma_k) \left(x+ u_k(x)\right)-x.
\end{equation*}
Observe that $u^{0}_{k+1}$ satisfies
the boundary condition
\begin{equation*}
u^{0}_{k+1}(x) = \sigma(\gamma_{k+1})x-x \quad\mbox{ for all } x \in \Gamma.
\end{equation*}
We use the P-nCG algorithm to compute a new local minimizer
$u_{k+1}$ ``near'' $u^{0}_{k+1}$, which satisfies the new
boundary condition $u_{k+1}(x) = \sigma(\gamma_{k+1})x-x$ for
all $x \in \Gamma$. We call this process the \emph{shear
loading continuation}.

We terminate the shear loading continuation when we can no longer find an equilibrium $u_{k+1}(x)$ close to the initial guess $\sigma(\del \gamma_{k+1})\left(x+ u_k(x)\right)-x.$ In practice, this means that we stop when the nCG algorithm returns a deformation in which the cores of the dislocations have moved apart.  If the dislocations move apart at the $(k+1)^{th}$ step of the iteration, then we assume that the critical strain must be between $\gamma_k$ and $\gamma_{k+1}$.

For the force based approximation~\eqref{eqn:qcf_force}, we use a similar iteration to find the critical strain.
At each step in the iteration, we solve the boundary value problem
\begin{equation}\label{eqn: force based continuation iteration}
\begin{alignedat}{3}
\mathcal{F}^{qcf}(u_{k+1})(x) &= 0  &&\quad\mbox{ for all } x \in \pdomain \setminus \Gamma,  \\
u_{k+1} (x) &= \sigma(\gamma_{k+1})x-x &&\quad\mbox{ for all }
x \in \Gamma.
\end{alignedat}
\end{equation}
We solve equation~\eqref{eqn: force based continuation iteration} using the
ghost force correction iteration defined in equation~\eqref{eqn: gfc iteration}.
We use the same initial guess
\begin{equation*}
u^{0}_{k+1}(x)  := \sigma(\del \gamma_{k+1})\left(x+  u_k(x)\right)-x
\end{equation*}
as in the energy-based case, and we solve each iteration of ghost force correction using the P-nCG algorithm.

We presented computational results for the deformation of a
Lennard-Jones chain under tension in~\cite{dobsonluskin08} that
demonstrate the necessity of using a sufficiently small parameter step
size to ensure that the computed solution remains in the domain of
convergence of the {\it ghost force correction iteration} (defined
below as QCE-QCF) method. These results exhibit fracture before the
actual load limit if the parameter step size is too large. We thus
conclude that the shear strain $\del \gamma_{k+1}$ must be
sufficiently small to ensure that our computed solution remains in the
domain of convergence of the QCE-QCF method since it would otherwise
predict instability for the Lomer dipole before that predicted by the
QCF method.

We perform the shear loading iteration for each of the quasicontinuum
methods using three different atomistic regions, as depicted in
Figure~\ref{fig: dipole atomistic regions}. The atomistic region (3)
(resp.,\ (4),\ (5)) is a box containing the dipole such that there are
three (resp., four, five) rows of atoms between the continuum region
and each of the atoms in the pentagons surrounding the cores of the
dislocations, or more precisely, the atomistic region ($k$) is the box
given in $\<\, \nu_1,\,\nu_2\,\>$ coordinates by
$[32-k,43+k]\times[30-k,30+k]$ for $k=3,4,5$.  Throughout the
remainder of this section, QNL-QCF (resp. QCE-QCF) will refer to the
QCF method implemented using the ghost force correction iteration
preconditioned by the QNL (resp. QCE) energy~\cite{qcf.iterative}. We
will let QNL-QCF(k) (resp.  QCE-QCF(k), QNL(k), QCE(k)) refer to the
QNL-QCF (resp.  QCE-QCF, QNL, QCE) method with atomistic region (k)
for k one of $3,\ 4, \text{ or } 5.$

We note that the atomistic regions studied in the benchmark
tests given in~\cite{Miller:2008} extended to the lateral
boundaries of the computational domain $\pdomain,$ or,
following the notation used in Figure~\ref{dipole} these tests
study the case $2D_1+w=L_1.$  Thus, the benchmark
tests~\cite{Miller:2008} do not present the capture of a Lomer
dipole in an atomistic region fully surrounded by a continuuum
region and do not study the effect of coupling error in the
dipole plane $\hkl(1 1 0).$

 In Table~\ref{tab: critical strain for
dipole}, we give the critical strains for each method in
decreasing order. The columns $\gamma^-$ and $\gamma^+$ are the
shear strains at the last step before the dislocations moved
apart, and at the first step for which the dislocations moved,
respectively. The column labeled ``Error" displays the relative
error from the critical strain of the atomistic model,
calculated by the formula
\begin{equation}\label{error}
 \text{Error} = \frac{\bar{\gamma}_{at} - \bar{\gamma}_{qc}}{\bar{\gamma}_{at}},
\end{equation}
where $\bar{\gamma}_{at}$ is the mean of $\gamma^-$ and $\gamma^+$ for
the atomistic method, and $\bar{\gamma}_{qc}$ is the mean of
$\gamma^-$ and $\gamma^+$ for the quasicontinuum method.  We wish to
stress that if $\gamma^+$ for some method is equal to $\gamma^-$ for
another method, then no strong statement can be made comparing
critical strains of the two methods.  For example, our data do not
show that QNL(5) has a significantly higher critical strain than
QCE-QCF(5).  Rather, our data suggest only that the critical strain of
the QCE-QCF(5) method is in the interval $[0.038172 ,0.038177]$ and
that the critical strain of the QNL(5) method is in $[0.038177,
0.038181]$.  Thus, the critical strains of the two methods are very
close.

The reader will observe that
the critical strain of the fully atomistic model is the highest.
The critical strain predicted by the QCE method is the lowest, and is the
farthest from the critical strain of the atomistic model.
This is not surprising in light of the results
on the stability of one dimensional chains obtained in ~\cite{doblusort:qce.stab}.
It is surprising, however, that the QNL method predicts the
critical strain of the atomistic model more accurately than
the QNL-QCF force-based method,
since we expected
that the stability of the QNL-QCF method would be determined by the QNL preconditioner~\cite{qcf.iterative}.
It also somewhat surprising that the QCE-QCF iteration predicted
the critical strain as accurately as the QNL-QCF iteration.
The analysis given in ~\cite{qcf.iterative} suggests that
the QCE-QCF iteration should lose stability at a lower strain than
the QNL-QCF iteration.

\begin{table}
\caption{The critical strain for the Lomer dislocation dipole
under shear.  The Error is given by $\left(\bar{\gamma}_{at} -
\bar{\gamma}_{qc}\right)/{\bar{\gamma}_{at}}$ in percent which
is defined precisely in the paragraph surrounding
\eqref{error}.} \label{tab: critical strain for dipole}
\begin{tabular}{l@{\qquad }c@{\quad}cr}
  Method & \multicolumn{2}{l}{Critical Strain} &  Error   \\
  \hline
  &$\gamma^{-}$ & $\gamma^{+}$ &  \\
Atomistic &   0.038190  &   0.038194  &     0.000 \% \\
   QNL(5) &   0.038177  &   0.038181  &     0.034 \% \\
 QCE-QCF(5) &   0.038172  &   0.038177  &     0.046  \% \\
 QNL-QCF(5) &   0.038172  &   0.038177  &     0.046 \% \\
   QNL(4) &   0.038164  &   0.038168  &     0.069 \% \\
 QCE-QCF(4) &   0.038155  &   0.038159  &     0.092 \% \\
 QNL-QCF(4) &   0.038155  &   0.038159  &     0.092 \% \\
   QNL(3) &   0.038111  &   0.038116  &     0.206 \% \\
QCE-QCF(3) &   0.038081  &   0.038085  &     0.286 \% \\
 QNL-QCF(3) &   0.038081  &   0.038085  &     0.286 \% \\
   QCE(3) &   0.037750  &   0.037775  &     1.125  \% \\
\end{tabular}
\end{table}

\begin{figure}\centering
 \subfigure{
  \includegraphics[trim=13cm 1cm 13cm .5cm, clip, width = 11.5 cm]{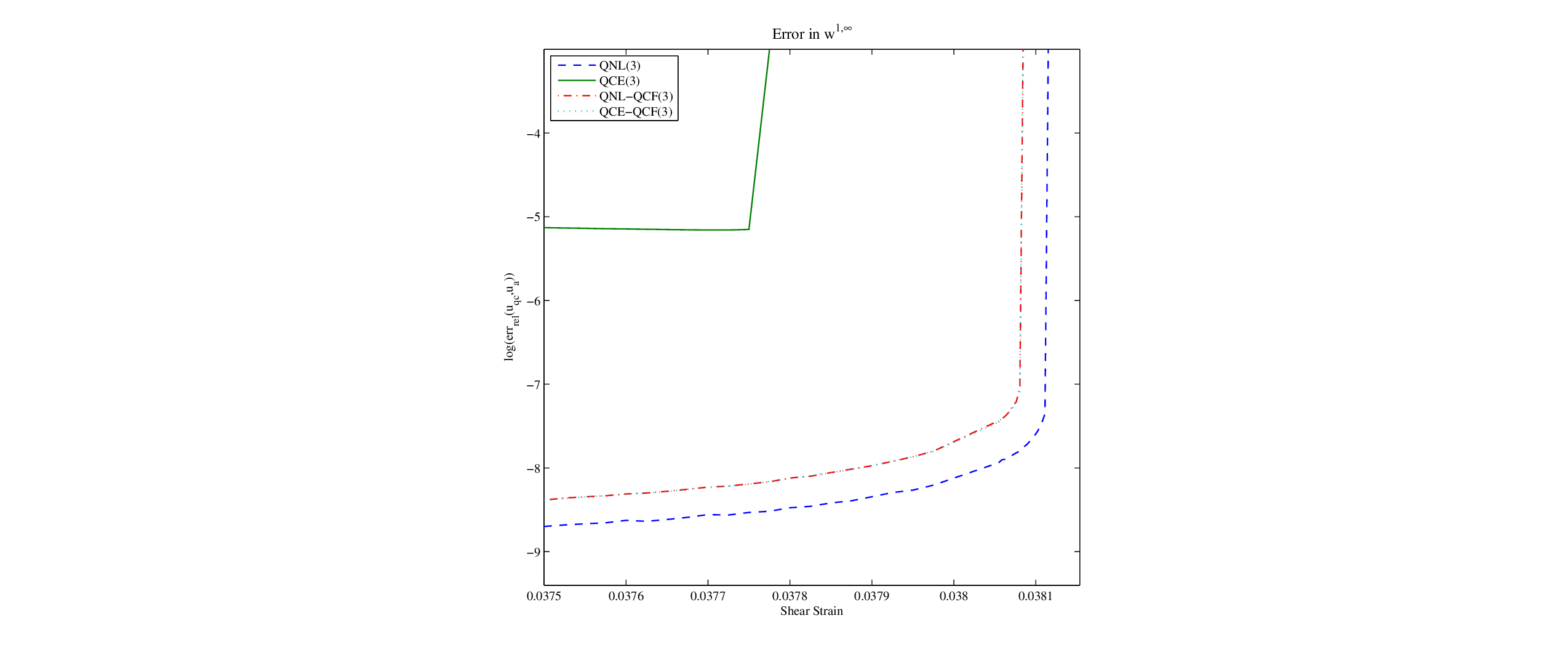}
  }
  \subfigure{
  \includegraphics[trim=13cm 2cm 13cm .5cm, clip, width = 11.5 cm]{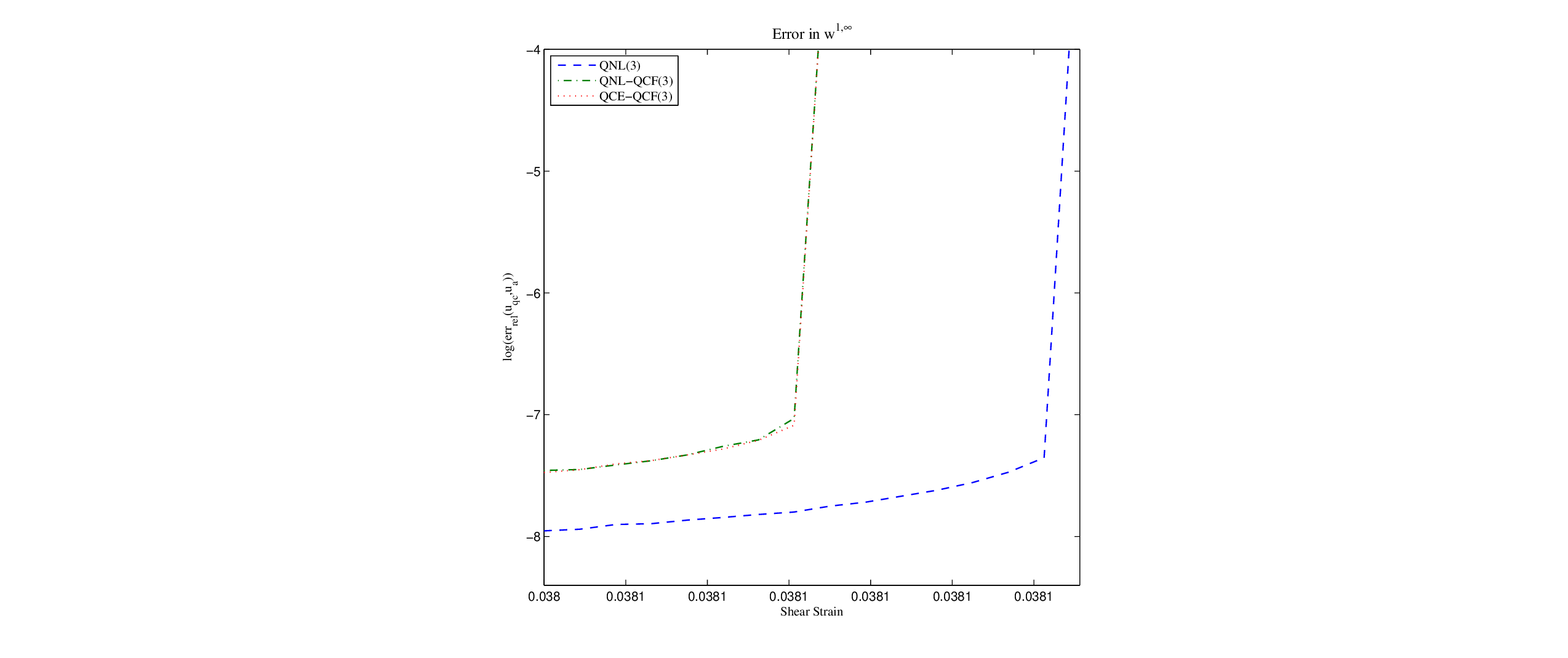}
  }
  \caption{The relative error $\rm{err_{rel}}(u_{qc}, u_a)=||u_{qc} - u_a||_{w^{1,\infty}}/||u_{a} ||_{w^{1,\infty}}$
  in $w^{1,\infty}$ for the Lomer dislocation dipole.}
  \label{fig: error for dipole}
\end{figure}

Figure~\ref{fig: error for dipole} consists of graphs showing the
relative $w^{1,\infty}$ error versus the shear strain for various quasicontinuum
methods. The
$w^{1,\infty}$ norm is defined by
\begin{equation*}
||u||_{w^{1,\infty}} := \max_{x\in\pdomain} \max_{b \in \mathcal{N}_x} \frac{|D_b u(x)|}{|b|},
\end{equation*}
and we define the \emph{relative error} in $w^{1,\infty}$ by
\begin{equation}\label{rel}
\rm{err_{rel}}(u_{qc}, u_a) := \frac{||u_{qc} - u_a||_{w^{1,\infty}}}{||u_{a}||_{w^{1,\infty}}}.
\end{equation}
Here $u_{qc}$ denotes a solution of one of the quasicontinuum models, and $u_a$ denotes
a solution of the atomistic model with the same boundary conditions.

The reader will observe that the error of the QCE method is the
greatest.  This is largely due to the presence of ghost
forces~\cite{Dobson:2008c}. The error of the QNL method is the least.
Again, this is surprising since various analyses ~\cite{Dobson:2008b,
  ortner:qnl1d, qcf.iterative, dobs-qcf2, MakrOrtSul} suggest that the
force based method should be more accurate than the QNL
method. However, we note that a similar effect was observed in
one-dimensional numerical experiments in \cite{MakrOrtSul}: while for
{large atomistic regions the QCF method was considerably more accurate than
the QNL method, for smaller atomistic regions the QNL method was clearly more
accurate.  The atomistic regions studied here in our numerical experiments for the
Lomer dislocation dipole are very small.}

\subsection{Tensile Loading}

In our next experiment, we examine the accuracy of various
quasicontinuum methods for the computation of the resistance
to slip under an applied tension.
By a tension, we mean a homogeneous deformation of the
lattice $\tri$ which takes the form
\begin{equation*}
  x \mapsto \tau(\gamma)x\quad\text{where}\quad\tau(\gamma):=\left (
    \begin{array}{cc}
     1+\gamma\quad & 0\\
      0 & 1
    \end{array} \right )
\end{equation*}
when expressed in the basis of the coordinate vectors $V_1$ and
$V_2$. We will let $\tau(\gamma)$ denote such a tension, and we
will call $\gamma$ the \emph{tensile strain}. We call a tension
\emph{positive} if $\gamma > 0$ and \emph{negative} if $\gamma
< 0$.

In our numerical experiments, we use the same pair potential
and lattice spacing as for the Lomer dipole problem.
We let $\pdomain$ be the rectangular domain consisting of $120$
columns of atoms parallel to $V_2$ with each column containing
$15$ atoms defined in the $\<\,\nu_1,\,\nu_2\,\>$ coordinates
by
\begin{equation*}
\begin{split}
\pdomain&=\left\{\< k,\,l\>: 0\le k\le 59,\ 0\le l\le 14\right\}\\
&\qquad\qquad\cup
\left\{\<k+1/2,\,l+1/2\>: 0\le k\le 59,\ 0\le l\le 14\right\}.
\end{split}
\end{equation*}
We will now let $\Gamma$ consist of four columns of atoms parallel to
$V_2$ at the left edge of $\pdomain$, and four columns at the right
edge of $\pdomain,$ that is,
\begin{equation*}
\begin{split}
\Gamma&=\left\{\<k,\,l\>: k=0,\,1,\,58,\,59,\ 0\le l\le 14\right\}\\
&\qquad\qquad\cup
\left\{\<k+1/2,\,l+1/2\>: k=0,\,1,\,58,\,59,\ 0\le l\le 14\right\}.
\end{split}
\end{equation*}
This corresponds to choosing Dirichlet boundary conditions on
the sides of $\omega$ parallel to $V_2$, and free boundary
conditions on the sides of $\omega$ parallel to $V_1$.

To find the critical strain for the crystal under tension,
we perform an iteration similar to the shear loading iteration.
Let $\gamma_k$ be an increasing sequence of tensile strains with
$\gamma_0 = 0$.  We let $u_0$ be the displacement field of the undeformed
reference lattice $\pdomain$
described above. Then we let $u_{k+1}$ solve the
problem
\begin{equation} \label{eqn: tension continuation problem}
\begin{alignedat}{3}
\mathcal{F}(u_{k+1})(x) &= 0  &&\quad\mbox{ for all } x \in \pdomain \setminus \Gamma,  \\
u_{k+1} (x) &= \tau(\gamma_{k+1})x -x&&\quad\mbox{ for all } x
\in \Gamma.
\end{alignedat}
\end{equation}
Here $\mathcal{F}$ is the gradient of one of the energies~\eqref{eqn: fully 2d energy},\ \eqref{eqn: qce energy},
\ \eqref{eqn: qnl energy}, or else the force-based method~\eqref{eqn:qcf_force}.
We solve equation~\eqref{eqn: tension continuation problem}
using the same methods discussed above for the shear loading iteration, except that
in this case we start each step with the initial guess
\begin{equation*}
u_{k+1}^{0}(x) := \tau \left (\smfrac{1 + \gamma_{k+1}}{1 + \gamma_k} \right ) \left(x+ u_k(x)\right)-x.
\end{equation*}
We call this iteration the \emph{tensile loading continuation}.
We stop the continuation when the P-nCG algorithm (or a ghost
force correction iteration) returns a deformation in which a
slip has occurred.

\begin{figure}
 \subfigure[The natural slip plane. We choose an atomistic region whose boundary
	    is the plane highlighted in grey. The mesh is the Delaunay triangulation \cite{CompGeom} of
	    the reference lattice.]{\includegraphics[width = 11.5 cm]{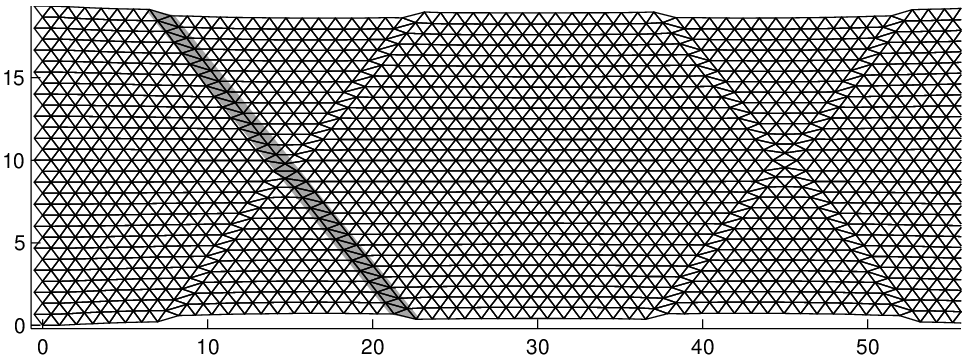}\label{fig: tension slip}}
 \subfigure[The atomistic region $\atomistic$ is shaded grey. All
 atoms in $\atomistic$ are marked with a blue dot.]{\includegraphics[width = 11.5 cm]{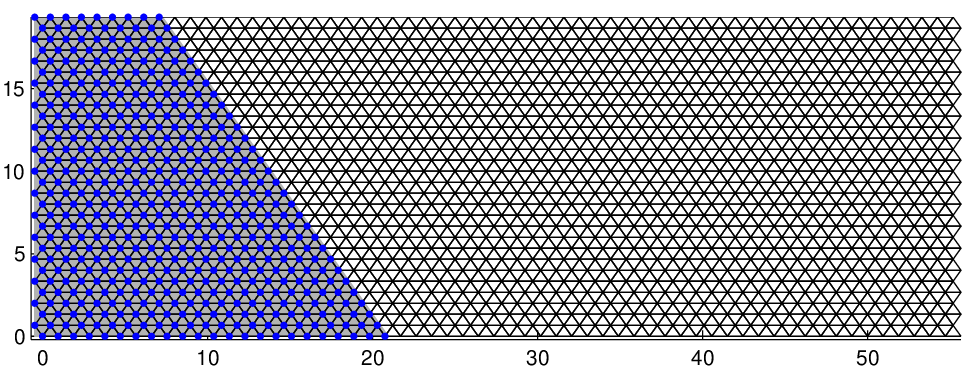}
 \label{fig: tension atomistic region}}
\caption{The atomistic ($\atomistic$) and continuum ($\mathcal{C}$)
  regions for the tensile loading experiment.}\label{slipping}
\end{figure}

Our numerical experiments are designed to test how well the boundary
between the atomistic and continuum regions resists slip under tensile
loading. When we perform the tensile loading iteration for the fully
atomistic model, we found that a slip tends to occur along the slip
plane indicated in Figure~\ref{fig: tension slip} when a critical
tensile load is reached.  Figure~\ref{fig: tension slip} shows the
configuration of the crystal immediately after a typical slip has
occurred.  We note that the slip allows the crystal to accommodate an
increased tensile strain with a decrease in the energy.  We then chose
an atomistic region $\atomistic$ whose boundary is along that line,
as depicted in Figure~\ref{fig: tension atomistic region}.

We observe that the critical strains of the quasicontinuum models were
lower than the critical strain of the atomistic model.  Our results
are summarized in Table~\ref{tab: critical strain for tension}.  In
Table~\ref{tab: critical strain for tension}, the columns $\gamma^-$
and $\gamma^+$ are the tensile strains at the last step before slip
is observed and at the first step for which slip is observed,
respectively. The column labeled ``Error" is the percent error from
the critical strain of the atomistic model.

\begin{table}
\caption{The critical strain under tension. The Error is given
    by $\left(\bar{\gamma}_{at} -
      \bar{\gamma}_{qc}\right)/{\bar{\gamma}_{at}}$ in percent, which
    is defined in
    \eqref{error}.} \label{tab: critical strain for tension}
\begin{tabular}{l@{\qquad }c@{\quad}cr}
  Method & \multicolumn{2}{l}{Critical Strain} &  Error   \\
  \hline
  &$\gamma^{-}$ & $\gamma^{+}$ &  \\
Atomistic &   0.081635  &   0.081649  &     0.000 \% \\
      QNL &   0.081593  &   0.081607  &     0.052 \% \\
  QCE-QCF &   0.081242  &   0.081270  &     0.473  \% \\
  QNL-QCF &   0.081101  &   0.081130  &     0.645 \% \\
      QCE &   0.063900  &   0.064200  &    21.548 \% \\
\end{tabular}
\end{table}

We observe that the critical strain of the QCE model is lower than the
critical strain of any of the other models. This is in agreement with
the one dimensional stability results established in
~\cite{doblusort:qce.stab}. However, we were surprised to find that
the critical strain predicted by the QNL method is higher than the
critical strain predicted by the QNL-QCF force based method, since we
expected that the stability of the QNL-QCF method would be determined
by the QNL preconditioner~\cite{qcf.iterative}.  We were also
surprised to observe that the critical strain of the QNL-QCF method
is less than the critical strain of QCE-QCF. The one-dimensional
analyses in ~\cite{doblusort:qcf.stab,qcf.iterative} suggest that the
force based methods should have a comparable critical strain as QNL,
and that the QCE-QCF iteration should lose stability at a lower strain
than the QNL-QCF iteration.

\begin{figure}\centering
 \subfigure{
  \includegraphics[trim=5cm .5cm 19.6cm 0cm, clip, width = 11.5 cm]{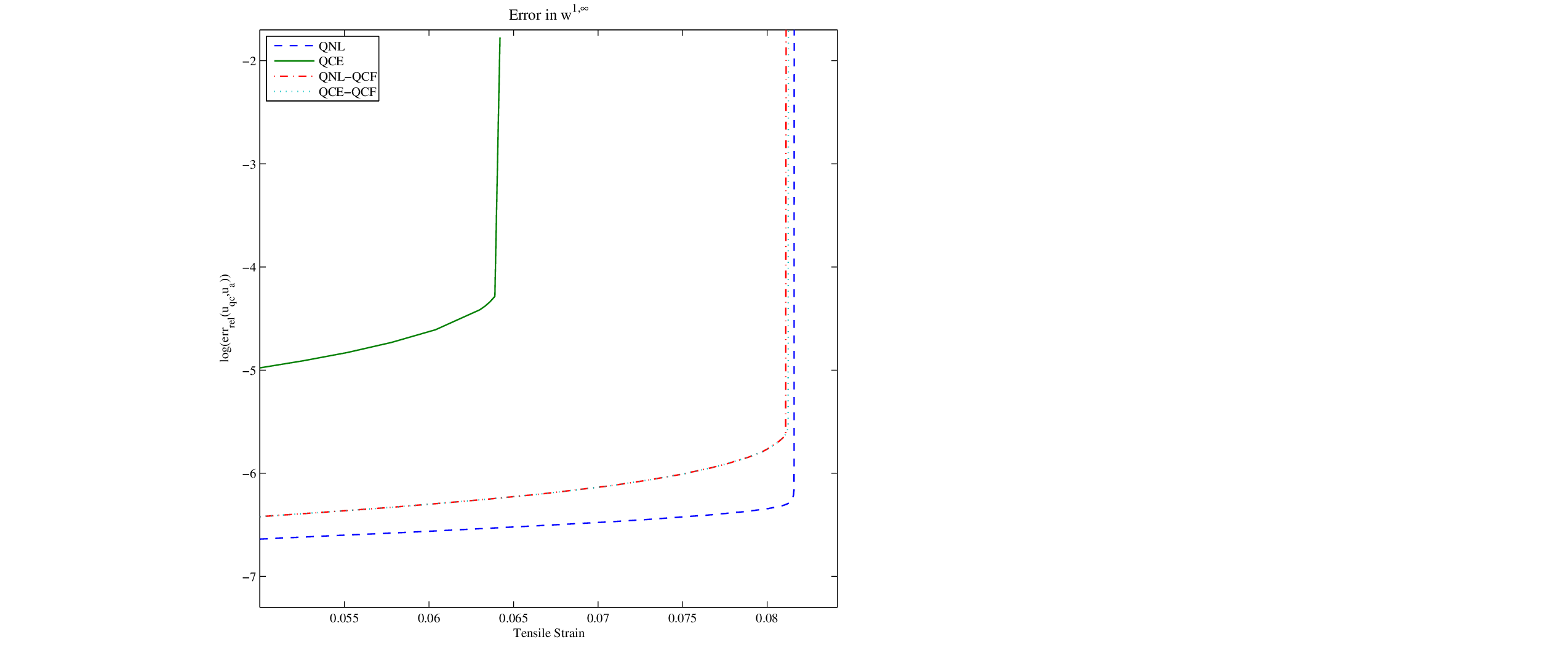}
  }
  \subfigure{
  \includegraphics[trim=13cm 1cm 13cm .5cm, clip, width = 11.5 cm]{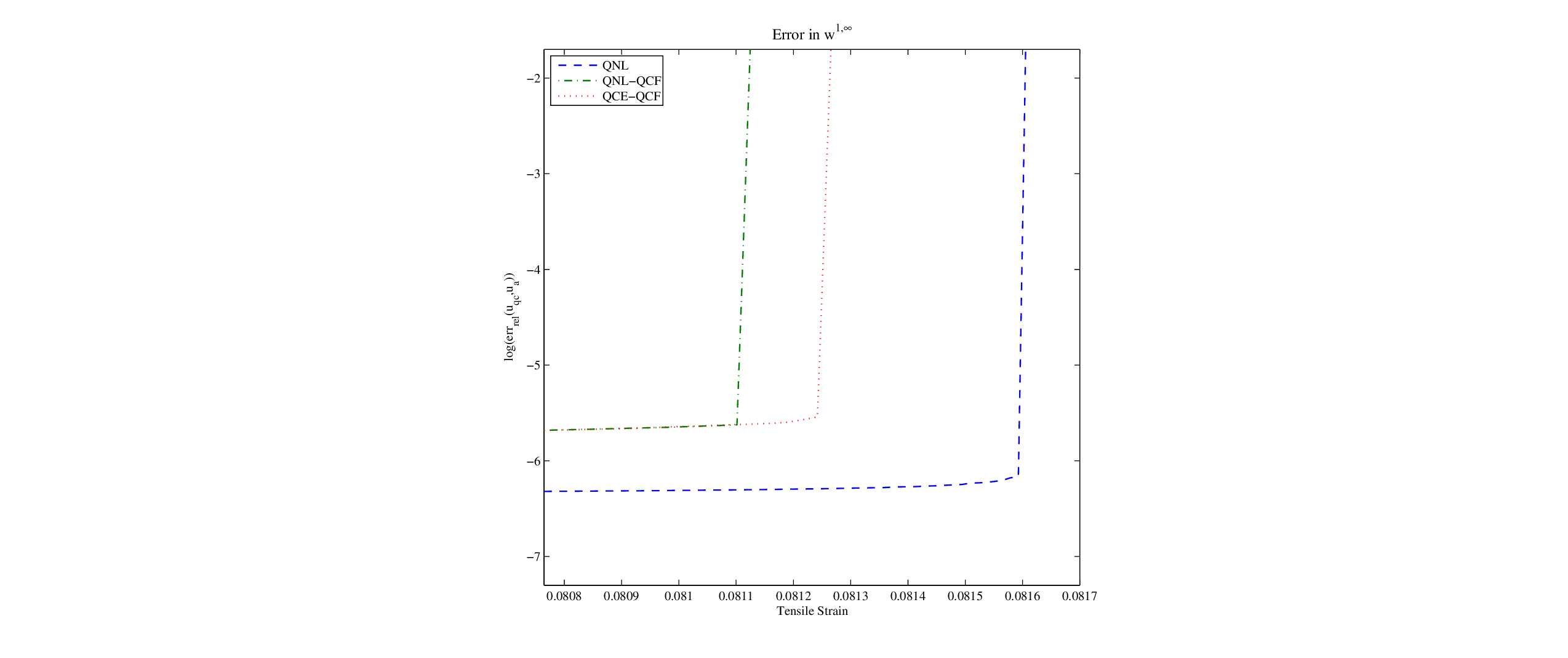}
  }
  \caption{The relative error $\rm{err_{rel}}(u_{qc}, u_a)=||u_{qc} - u_a||_{w^{1,\infty}}/||u_{a} ||_{w^{1,\infty}}$
  in $w^{1,\infty}$ for tensile loading.}
\label{fig: error for tension}
\end{figure}

Figure~\ref{fig: error for tension} consists of graphs showing the
relative $w^{1,\infty}$ error versus the tensile strain for various
quasicontinuum methods. We observe that the QCE model had the greatest
error. This is primarily the result of the ghost forces which arise in
the QCE model. The QNL model had the lowest error. Again, this is
somewhat surprising in light of the one dimensional
theory~\cite{Dobson:2008b, ortner:qnl1d, qcf.iterative, dobs-qcf2,
  MakrOrtSul}, which predicts that the force based methods should have
lower error than the QNL method.

\section{Conclusion}
Materials scientists and engineers typically attempt to verify
their multiphysics methods by benchmark tests.  However,
definitive conclusions from these benchmark tests are often not
possible since they combine many sources of error (modeling
error, coupling error, solution error, boundary condition
error, etc.).  Further supporting theory is needed to select a
set of test problems that thoroughly samples the solution
space.

We are developing an improved theoretical basis for
benchmarking atomistic-to-continuum coupling methods based on
the multiscale numerical analysis developed by us and others.
Our goal is to be able to reliably predict the accuracy of
atomistic-to-continuum coupling methods for general
deformations and loads from numerical experiments for a small
set of mechanics problems. This set of mechanics problems
should sample the fundamental modes of material instability
such as dislocation formation, slip, and fracture.

Our discussion of the benchmark tests presented in this paper
give many examples of the predictive success of the theoretical
analysis we have developed during the past few years,
but we also describe several cases where our theoretical
analysis seems to predict a different outcome than our
computational experiments.  This discrepancy between theory
and computational experiment occurs when our theoretical
analysis does not adequately
model the computational problem and is motivation to develop
more general theoretical analysis.

\section{Appendix: A Preconditioned Conjugate Gradient Algorithm}
\label{sec:pcg}
In this appendix, we describe the preconditioned nonlinear conjugate
gradient optimization algorithm (P-nCG) and the corresponding
linesearch method, which we used in the numerical experiments in this
paper.

Let $\E : \R^N \to \R\cup\{+\infty\}$ be continuously differentiable
in $\{u : \E(u) < \infty \}$. If $\E$ is an energy of the type
discussed above, then the standard nonlinear conjugate gradient
method~\cite[Sec. 5.2]{NocWri} is convergent, but very inefficient,
due to the poor conditioning of the Hessian matrix at local minima.

Let ${\rm P} \in \R^{N \times N}$ be a symmetric positive
definite matrix, e.g., a discrete Laplacian or a modification
thereof. A simple but considerably more efficient method is
obtained if all inner products in the conjugate gradient
algorithm are replaced by a ${\rm P}$-inner product $(u,
v)_{\rm P} = u^T {\rm P} v$, and all gradients $\del \E(u)$ by
${\rm P}$-gradients $\del_{\rm P} \E(u) = {\rm P}^{-1} \del
\E(u)$. The ${\rm P}$-gradient $\del_{\rm P}\E(u)$
``represents'' the gradient $\del \E(u)$ in the $P$-inner
product $(v,w)_{\rm P}:=v^T{\rm P}w$, since
\begin{displaymath}
  (\del_{\rm P} \E(u), w)_{\rm P} = ({\rm P}^{-1} \del \E(u))^T {\rm P} w = \del \E(u)^T
  w = (\del \E(u), w).
\end{displaymath}
In practice, we allow a new preconditioner to be computed at
each step. A basic preconditioned conjugate gradient algorithm
of Polak-Ribi{\`e}re type can be
described as follows: \\[2mm]
{
.}
\qquad
\begin{minipage}{0.8\textwidth}
\begin{enumerate}
\item[(0)] Input: $u_0 \in \R^N$;
\item[(1)] Evaluate ${\rm P}_0$;\ \ $g_0 = {\rm P}_0^{-1}
    \del \E(u_0)$;\ \ $s_0 = 0$;
\item[(2)] For $n = 1, 2, \dots$ do:
\item[(3)] \qquad Evaluate ${\rm P}_n$;
\item[(4)] \qquad $g_n = {\rm P}_n^{-1} \del\E(u_n)$;
\item[(5)] \qquad $\beta_n = \max\{0, (g_n, g_n - g_{n-1})_{{\rm P}_n} / (g_{n-1},
  g_{n-1})_{{\rm P}_{n-1}} \}$;
\item[(6)] \qquad $s_n = - g_n + \beta_n s_{n-1}$;
\item[(7)] \qquad $\alpha_n \leftarrow {\rm LINESEARCH}$;
\item[(8)] \qquad $u_n = u_{n-1} + \alpha_n s_n$;
\end{enumerate}
\vspace{2mm}
\end{minipage}

In the following we specify further crucial details of our
implementation:
\begin{enumerate}
\item {\it Preconditioner: } The choice of preconditioner has the
  biggest influence on the efficiency of the optimization. We project
  all atoms onto a single plane and triangulate the resulting set of
  nodes. {On this triangulation, we assemble the standard P1-finite
  element stiffness matrix, ${\rm K}_n$, discretizing the negative
  Laplace operator. The preconditioner ${\rm P}_n$ is obtained by
  imposing homogeneous Dirichlet boundary conditions on the clamped
  nodes.}

\item {\it LINESEARCH: } Our linesearch is implemented following
  Algorithms 3.5 and 3.6 in \cite{NocWri} closely. We use cubic
  interpolation from function and gradient values in the
  ``Interpolate'' step of Algorithm 3.6. We guarantee the strong Wolfe
  conditions, \cite[Eq. (3.7)]{NocWri} with constants $c_1 = 10^{-4}$,
  $c_2 = 1/2$.

\item {\it Initial guess for $\alpha_n$: } If the initial guess
  $\alpha_n^{(0)}$ for the new steplength $\alpha_n$, which is passed
  to the LINESEARCH routine, is chosen well, then actual linesearch
  can be mostly avoided, which can result in considerable performance
  gains. Following \cite[Eq. (3.60)]{NocWri}, and extensive
  experimentation with alternative options, we choose $\alpha_n^{(0)}
  = 2 (\E(u_{n-1}) - \E(u_{n-2})) / (g_{n}, s_{n})_{P_n}$.

\item {\it Termination Criteria: } We terminate the
    iteration successfully if the following condition is
    satisfied:
  \begin{align*}
    & \Big( \|u_n - u_{n-1}\|_\infty \leq {\rm TOL}_u^{\infty} \quad
    \text{or} \quad \|u_n - u_{n-1}\|_{{\rm P}_n} \leq {\rm TOL}_u^P \Big) \\
    {\rm and} \quad & \Big( \|\del \E(u_n)\|_\infty \leq {\rm
      TOL}_g^\infty
    \quad \text{or} \quad \|g_n\|_{{\rm P}_n} \leq {\rm TOL}_{g}^{P} \Big) \\
    \text{and} \quad &\Big( \E(u_{n-1}) - \E(u_n) \leq {\rm
      TOL}_\E \Big),
  \end{align*}
  where $\|\cdot\|_\infty$ denotes the $\ell^\infty$-norm,
  and $\|\cdot\|_{\rm P}$ denotes norm associated with the
  ${\rm P}$-inner product. The tolerance parameters are adjusted
  for each problem. Typical choices are ${\rm TOL}_u^\infty
  = {\rm TOL}_u^P = 10^{-5},\ {\rm TOL}_g^\infty = {\rm
    TOL}_g^P = 10^{-4},\ {\rm TOL}_\E = 10^{-4}$.

  We terminate the iteration unsuccessfully if a maximum number of
  iterations is reached, or if the LINESEARCH routine is unable to
  make any progress. Note also that, since P-nCG is a descent method,
  we have $\E(u_{n-1}) - \E(u_n) = |\E(u_{n-1}) - \E(u_n)|$.

\item {\it Robustness Checks: } In addition, our algorithm
    uses various minor modifications to increase its
    robustness. Since these do not significantly affect its
    performance, we have decided to not give any details.
\end{enumerate}

The algorithm described above is both efficient and robust for most of the problems we consider. However, we wish to stress certain difficulties that arise in the presence of ``meta-stable'' states and particularly shallow local minimizers, which are difficult to distinguish numerically. In our experience, dislocations fall precisely into this category.

On some occasions, our algorithm would fail when a particularly
low tolerance setting was used. The reason for the failure is
usually that the directional derivative along the search
direction is non-negative (up to numerical precision) and hence
the linesearch fails or stagnates. Replacing the conjugate
direction by a steepest descent direction (one of our
robustness checks) resolves this problem only partially.  We
were able to overcome these difficulties mostly by tweaking the
various optimization parameters.

\section{Acknowledgements}
We wish to thank Ellad Tadmor for helpful discussions.

\end{document}